\documentclass[11pt]{amsart}
\usepackage{amssymb,amsmath,amsthm,enumitem,xcolor}
\IfFileExists{lmodern.sty}{\usepackage[T1]{fontenc}\usepackage{lmodern}}{}
\usepackage[colorlinks=true,linkcolor=blue,citecolor=blue,urlcolor=blue]{hyperref}
\usepackage[margin=1.15in]{geometry}
\numberwithin{equation}{section}
\newtheorem{thm}{Theorem}[section]
\newtheorem{lem}[thm]{Lemma}
\newtheorem{prop}[thm]{Proposition}
\newtheorem{cor}[thm]{Corollary}
\theoremstyle{definition}
\newtheorem{defn}[thm]{Definition}
\newtheorem{rem}[thm]{Remark}
\newcommand{\Q}{\mathbf Q}
\newcommand{\Z}{\mathbf Z}
\newcommand{\F}{\mathbf F}
\newcommand{\A}{\mathbf A}
\newcommand{\cO}{\mathcal O}
\newcommand{\cA}{\mathcal A}
\newcommand{\cD}{\mathcal{D}}
\newcommand{\cC}{\mathcal C}
\newcommand{\GL}{\mathrm{GL}}
\newcommand{\PSL}{\mathrm{PSL}}
\newcommand{\PGL}{\mathrm{PGL}}

\newcommand{\eps}{\varepsilon}
\newcommand{\Qbar}{\overline{\Q}}
\newcommand{\Rb}{\overline R}
\newcommand{\Sb}{\overline S}
\newcommand{\gt}{\tilde g}

\begin{document}
\title[How to facet a gemstone]{How to facet a gemstone: from potential modularity to modularity lifting}
\author{Luis Dieulefait}
\date{2 December 2026}
\begin{abstract}
Let $F$ be a totally real field and $p\ge7$ a prime, possibly ramified in $F$. We prove that a potentially modular $p$-adic representation of $G_F$ which is potentially crystalline with distinct Hodge--Tate weights above $p$, and whose reduction is modular and remains absolutely irreducible with non-dihedral projective image on $G_{F(\zeta_p)}$, is modular. This removes the hypotheses of \cite{MF} on the Serre weights, namely the multiplicity-free hypothesis and the exclusion of Steinberg weights, and its hypothesis that $p$ is unramified in $F$. The same method applies to potentially semistable representations. The proof combines the identity of cycles of \cite{MF} with the modularity of lifts whose local components are mixed along an extension in which $p$ splits. These lifts exist by a theorem of Emerton, Gee, Pan and Zhu, they are potentially automorphic by a theorem of Calegari, Emerton and Gee, and their modularity is proved by tensoring with lifts of an auxiliary residual representation of a different weight: the tensor products lie on a common component of a local deformation ring for $\GL_4$, a lifting theorem of Thorne applies, and the factors of an automorphic tensor product are modular by \cite{Die26}.
\end{abstract}
\maketitle

\section{Introduction}
\subsection{The results}
This paper is a sequel to \cite{MF}, and we use its notation. Let $F$ be a totally real field and $p$ a prime, which may ramify in $F$. For a place $v\mid p$ and an integer $k\ge2$, $W_{k,v}$ is the semisimplified reduction of $\bigotimes_{\kappa\colon F_v\hookrightarrow\Qbar_p}\mathrm{Sym}^{k-2}$, where $\GL_2(\cO_{F_v})$ acts on the factor indexed by $\kappa$ through $\kappa$; it factors through $\GL_2(k_v)$, and for unramified $F_v/\Q_p$ it is the semisimplification of the coefficient representation in \cite[\S1.1]{MF}; for $\bar r\colon G_{F_v}\to\GL_2(\overline\F_p)$, $W(\bar r)$ is its set of Serre weights; and a Serre weight is Steinberg if it is a twist of $\bigotimes_\tau\mathrm{Sym}^{p-1}$.

\begin{thm}\label{thm:main}
Let $p\ge7$ be a prime, $F$ a totally real field, and $k\ge2$. Let $\rho\colon G_F\to\GL_2(\Qbar_p)$ be a continuous representation, unramified outside finitely many places, with reduction $\Rb$. Assume:
\begin{enumerate}[label=\textup{(\arabic*)}]
\item for every $v\mid p$, $\rho|_{G_{F_v}}$ is crystalline with Hodge--Tate weights $0$ and $k-1$ at every embedding;
\item $\rho$ is potentially modular: $\rho|_{G_{F'}}$ comes from a Hilbert modular form for some finite Galois totally real extension $F'/F$;
\item $\Rb$ is modular, and $\Rb|_{G_{F(\zeta_p)}}$ is absolutely irreducible with non-dihedral projective image.
\end{enumerate}
Then $\rho$ is modular: it comes from a Hilbert modular eigenform over $F$ of parallel weight $k$ and level prime to $p$.
\end{thm}

In \cite[Theorem 1.2]{MF} there was a fourth hypothesis: every Serre weight of $\Rb|_{G_{F_v}}$ occurring in $W_{k,v}$ had to be non-Steinberg and to occur with multiplicity one. That condition restricted $k$ to certain ranges; here $k$ is arbitrary and there is no condition on the Serre weights. As in \cite{MF}, nothing is assumed at the places away from $p$.

The general statement is the following. For a regular Hodge type $\lambda_v$ and an inertial type $\tau_v$ at $v\mid p$, $\sigma(\lambda_v,\tau_v)$ is the locally algebraic type of \cite[\S1.2]{MF} and $\overline{\sigma(\lambda_v,\tau_v)}$ the semisimplified reduction of a lattice in it.

\begin{thm}\label{thm:general}
Let $p\ge7$ be a prime and $F$ a totally real field. Let $\rho\colon G_F\to\GL_2(\Qbar_p)$ be a continuous representation, unramified outside finitely many places, with reduction $\Rb$. Assume:
\begin{enumerate}[label=\textup{(\arabic*)}]
\item for every $v\mid p$, $\rho|_{G_{F_v}}$ is potentially crystalline with distinct Hodge--Tate weights at every embedding; let $\lambda_v$ and $\tau_v$ be its Hodge type and its inertial type;
\item $\rho$ is potentially modular;
\item $\Rb$ is modular, and $\Rb|_{G_{F(\zeta_p)}}$ is absolutely irreducible with non-dihedral projective image.
\end{enumerate}
Then $\rho$ is modular: it comes from a Hilbert modular eigenform over $F$, of the weight given by the $\lambda_v$, whose local component $\pi_v$ at each $v\mid p$ contains $\sigma(\tau_v)^\vee$.
\end{thm}

Theorem~\ref{thm:main} is the case of trivial inertial types and parallel weight. Theorem~\ref{thm:general} extends to potentially semistable representations (Theorem~\ref{thm:ss}). The weights in Theorem~\ref{thm:general} need not be parallel. For Hodge--Tate weights $0$ and $1$ the statement follows from the potentially Barsotti--Tate results; see \cite[\S5.4]{GK}.

\begin{rem}[On the hypotheses]
(a) Hypothesis (2) plays the same role as in \cite{MF}; see \cite[Remark 1.3]{MF}. It holds for members of compatible systems that are potentially modular at some prime.

(b) There is no hypothesis on the Serre weights of $\Rb|_{G_{F_v}}$. In \cite{MF} Steinberg weights were excluded, because the identity of cycles was proved there with tame types, which never contain a Steinberg weight, and because the cycle of a Steinberg weight can share components with cycles of other weights. The argument here does not require disjointness or coefficient one. Here only the fact that the cycle of each Serre weight is a local object is used, and for Steinberg weights this is obtained in Section~\ref{sec:St} from one family of wildly ramified types. In particular Theorem~\ref{thm:main} includes the weight $k=p+1$.

(c) The bound $p\ge7$ comes from the results of \cite{CEG} as they are stated. The lifting theorem that we use \cite[Theorem 7.1]{Tho12} requires an adequate residual image; for the four-dimensional representations to which we apply it this is proved in Lemma~\ref{lem:aux}, using the results on tensor products of \cite[\S2.2]{ADP}.

(d) There is no hypothesis on the ramification of $p$ in $F$. In \cite{MF} it was assumed that $p$ is unramified in $F$; as observed in \cite[Remark 6.5]{MF}, this was used for the combinatorics of Serre weights and to simplify the notation. The results of \cite{MF} used here hold without it; see \S\ref{ss:ramified}. When $p$ ramifies at $v$, several embeddings of $F_v$ induce the same embedding of $k_v$, so the same residual embedding can contribute several symmetric powers to $W_{k,v}$, which is harmless since no multiplicity condition is imposed.

(e) In hypothesis (3), a projective image isomorphic to the Klein four-group is regarded as dihedral. The hypothesis is used to make all lifts Lie-irreducible, as required in \cite{Die26}, and in the proof of adequacy.
\end{rem}

\subsection{The idea}
Recall the argument of \cite{MF}. After a soluble base change one is in a \emph{global setting} over $F$, where patching gives an identity between cycles on the special fibre of the local deformation rings above $p$:
\begin{equation}\label{eq:idF}
\sum_{C\in\cA}\ \prod_{v\mid p}Z(\overline C_v)=\prod_{v\mid p}s_v .
\end{equation}
Here $\cA$ is the set of \emph{automorphic components}, the families $C=(C_v)_v$ of irreducible components of the crystalline deformation rings that contain the local representations of an eigenform lifting $\Rb$; $Z(\overline C_v)$ is the cycle of the special fibre of $C_v$; and $s_v$ is a cycle depending only on $\Rb|_{G_{F_v}}$ and $k$, the sum of the cycles of the Serre weights in $W_{k,v}$ with their multiplicities. If $L/F$ is a totally real extension of degree $d$ in which the places above $p$ split completely, the same identity holds over $L$, and its right-hand side is the $d$-th tensor power of the right-hand side of \eqref{eq:idF}.

Take $L$ such that $\rho$ is modular over $L$. The left-hand side over $L$ contains the term of the constant family of the components of $\rho$. If we knew that \emph{every} family over $L$ obtained by choosing, sheet by sheet, an element of $\cA$ is automorphic over $L$, these $|\cA|^d$ families would already account for the whole right-hand side, by the $d$-th power of \eqref{eq:idF}; no room would be left, so the family of $\rho$ would be one of them, and the components of $\rho$ would be in $\cA$. Then $\rho$ is modular over $F$ by a standard patching argument. In \cite{MF} the multiplicity-free hypothesis made this unnecessary: there the relevant coefficient is $1$ and the constant families suffice.

So the problem is to prove that lifts of $\Rb|_{G_L}$ whose local components are \emph{mixed}, equal to those of one eigenform at some places above $p$ and to those of another at other places, are modular over $L$. Such lifts exist \cite{EGPZ} and are potentially automorphic \cite{CEG}. To prove their modularity we use an auxiliary residual representation $\Sb$, unrelated to $\Rb$, coming from an eigenform $g$ of a weight different from $k$. After a soluble extension $F'/F$, ramified above $p$, the restrictions of $\Rb$ and $\Sb$ to the decomposition groups above $p$ become trivial, so a local component $C$ of a lift of $\Rb$ and the local component $E$ of $g$ are components of deformation rings of the same residual representation. If $A$ lifts $\Rb$ and $B$ lifts $\Sb$, the local representation of $A\otimes B$ at a place where $A$ is on $C$ and $B$ on $E$, and at a place where $A$ is on $E$ and $B$ on $C$, lie on the \emph{same} component of the crystalline deformation ring for $\GL_4$: the tensor product does not see the order of the factors. We say that two such representations are \emph{potentially connected} (Definition~\ref{def:potconn}): they need not have the same reduction over the base field, but over a finite extension they lie on a common component. The two weight differences are distinct, so the tensor has four distinct Hodge--Tate weights. By Lemma~\ref{lem:aux} the residual tensor is absolutely irreducible with adequate image on $G_{F(\zeta_p)}$. Hence a lifting theorem on a fixed component \cite{Tho12} transfers automorphy from $A_0\otimes B_0$ to $A\otimes B$ whenever, place by place, the two pairs of components agree up to order. By \cite{Die26}, the automorphy of $A\otimes B$ gives the modularity of $A$ and $B$. Starting from the pair formed by an eigenform lifting $\Rb$ and $g$, an induction on the number of distinct components produces all mixed lifts.

\subsection{Antecedents of the tensor trick}
The use of tensor products to move a representation into a situation where a lifting theorem applies goes back to the tensor product trick of Harris \cite{Har09}, in the form developed by Barnet-Lamb, Gee, Geraghty and Taylor \cite{BLGGT}. There one tensors with a monomial representation, the induction of a character of a CM extension, chosen so that the local components of the tensor product above $p$ are under control; the original representation is recovered from the tensor product because the auxiliary factor is induced. The control of the local components is achieved in \cite{BLGGT} through potential diagonalizability, and the relation of lying on a common component of a local lifting ring is the relation $\sim$ of \cite{BLGGT}. The argument of this paper is inspired by these ideas, and the notion of potentially connected representations of \S\ref{ss:potconn} is a variant of that relation.

Two things are different here. First, the trick is applied to components that are not known to be potentially diagonalizable: nothing is used about the components of the two factors except that the tensor product does not depend on their order. Second, the auxiliary factor is not monomial. It is the representation of an eigenform without complex multiplication and with large residual image, so automorphic induction is not available; its place is taken, in one direction, by the automorphy of the Rankin--Selberg product of Ramakrishnan \cite{Ram00}, and in the other direction by its converse \cite{Die26}, which recovers the modularity of the factors from the automorphy of the tensor product.

\subsection{Organization} Section~\ref{sec:inputs} states the results we use. Section~\ref{sec:aux} constructs the auxiliary form. Section~\ref{sec:setup} fixes the fields and the local components, defines potentially connected representations and proves the lemma on tensor products. Section~\ref{sec:mixed} proves the modularity of mixed lifts. Section~\ref{sec:St} extends the identity of cycles of \cite{MF} to Steinberg weights. Section~\ref{sec:proof} proves Theorem~\ref{thm:main}, Section~\ref{sec:general} explains the changes for Theorem~\ref{thm:general}, and Section~\ref{sec:ss} those for potentially semistable representations.

\subsection{Conventions} As in \cite{MF}. The cyclotomic character $\eps$ has Hodge--Tate weight $-1$, and $\operatorname{Art}_K$ sends uniformizers to arithmetic Frobenius. In the automorphic correspondence, Frobenius elements are geometric. We use L-algebraic automorphic representations when their unramified Satake parameters are the eigenvalues of the associated Galois representation. If $\Pi$ is such a representation of $\GL_n$, its C-algebraic normalization is $\Omega=\Pi|\det|^{(n-1)/2}$, so that $r_p(\Omega)$ has these eigenvalues \cite[Theorem 2.1.1]{BLGGT}. A representation of $G_{F'}$, $F'$ totally real, is \emph{modular} if it comes from a Hilbert modular eigenform, up to a twist by a character of the form (finite order)$\cdot\eps^n$; it is \emph{Lie-irreducible} if it is irreducible on every open subgroup; a residual representation is \emph{dihedral} if it is induced from a character of a quadratic extension.

\section{Results used}\label{sec:inputs}
\begin{enumerate}[label=\textup{(I\arabic*)},leftmargin=3em]
\item \emph{The comparison of \cite{MF}.} The reduction to a global setting \cite[\S4.1]{MF}; the properties $(\dagger)$, the module $M_k$, the set $\cA$ of automorphic components and \cite[Lemma 3.6]{MF}: a family $C$ is in $\cA$ if and only if some modular lift with the properties $(\dagger)$ has its local points on $C$, and a lift with the properties $(\dagger)$ whose local points lie on a family of $\cA$ comes from an eigenform contributing to $M_k$; the cycle identity \cite[Lemma 3.7]{MF}; the independence of the classes of tuples \cite[Lemma 2.6]{MF}; the Steinberg ring \cite[Proposition 2.4]{MF}; and \cite[Lemma 4.1]{MF}, which rests on \cite[Theorem 2.1.16]{CEG}.
\item \emph{Lifts on prescribed components.} \cite[Theorem 1.1.1]{EGPZ}: let $F'$ be totally real and $\bar r\colon G_{F'}\to\GL_2(\F)$ totally odd with $\bar r|_{G_{F'(\zeta_p)}}$ absolutely irreducible, $p>5$. For every global deformation problem, consisting of a de Rham character lifting $\det\bar r$, a finite set of places, a non-empty union of irreducible components of the fixed-determinant lifting ring at each of these places away from $p$, and a non-empty union of irreducible components of a potentially semistable lifting ring of regular Hodge type at each place above $p$, the global deformation ring is finite over $\cO$ of positive rank; in particular lifts of that type exist. More precisely \cite[Definition 2.1.2]{EGPZ}: the set of places contains those above $p$ and those where $\bar r$ or the character is ramified, and lifts are unramified outside it; at $v\mid p$ one fixes a regular Hodge type $\lambda_v$ and a finite extension $E_v$ of the completion, and takes a union of components of the ring of lifts of determinant $\psi$ and Hodge type $\lambda_v$ that become semistable over $E_v$, which is assumed to be non-zero. There is no condition on the ramification of $p$ in $F'$ or on the local fields, and the Hodge types at different places are independent. The theorem does not assert that the lifts are potentially modular; for that we use (I3).
\item \emph{Potential automorphy from realizable components.} \cite[Theorem B, Corollary 4.2.12]{CEG}, for $n=2$: let $E_0$ be a CM field and $(s,\mu)$ a polarized, odd representation of $G_{E_0}$, ramified at finitely many primes and reasonable, such that the local representations at the places above $p$ are globally realizable. Then $(s,\mu)$ is automorphic over a finite Galois CM extension of $E_0$. Here \emph{reasonable} means: $p>6$, $\zeta_p\notin E_0$, and the reduction is polarizable, odd and irreducible on $G_{E_0(\zeta_p)}$ \cite[Definition 2.1.6]{CEG}. A component of a potentially semistable lifting ring over a local field $K$ is \emph{globally realizable} if it contains the local representation, at a place with completion $K$, of a member with reasonable reduction of an odd, regular, polarized, weakly irreducible compatible system over some CM field, with the same local multiplier and the same local residual representation \cite[Definition 2.1.9]{CEG}; the global residual representation of that system need not be related to $\bar s$. Components are irreducible components of the generic fibres of the rings of potentially semistable lifts of a fixed inertial type and a fixed regular Hodge type \cite[\S1.5]{CEG}. At a place above $p$ that splits in $E_0$, the polarized lifting ring is identified with the lifting ring for $\GL_n$, and global realizability does not depend on the local multiplier \cite[Remark 2.1.11]{CEG}. The proof of \cite[Corollary 4.2.12]{CEG} gives the extension linearly disjoint from any given finite extension of $E_0$. The theorem gives potential automorphy only; it does not give automorphy even when another member of the compatible system is known to be automorphic \cite[\S1.1]{CEG}.
\item \emph{Lifting on a fixed component.} We use the following sufficient form of \cite[Theorem 7.1]{Tho12} and \cite[Corollary 7.3]{Tho17}. Let $E_0$ be a CM field with $\zeta_p\notin E_0$ and $r\colon G_{E_0}\to\GL_n(\cO)$ with $r^c\simeq r^\vee\eps^{1-n}\mu$, where $\mu$ is a character of the Galois group of the maximal totally real subfield taking the same value at all complex conjugations. Assume that $r$ is ramified at finitely many places and that $\bar r$ is absolutely irreducible with $\bar r(G_{E_0(\zeta_p)})$ adequate. Assume that there is $r'$, with the same reduction and multiplier modulo $p$, attached to a regular algebraic, essentially conjugate self-dual, cuspidal automorphic pair $(\pi,\chi)$, where $\pi$ is a representation of $\GL_n(\A_{E_0})$ potentially unramified above $p$ and the reduction of $r_p(\chi)$ is $\bar\mu$, such that $r'\sim r$ at the places above $p$, and such that at each other place either $\pi$ and $r$ are unramified or $r'\leadsto r$. Then $r$ is automorphic. Here $r'\sim r$ above $p$ means that the two local representations are points of the same irreducible component of the generic fibre of a potentially crystalline lifting ring with fixed Hodge--Tate weights \cite[Definition 3.7]{Tho12}; away from $p$, $r'\sim r$ means the same for the unrestricted lifting ring, and $r'\leadsto r$ means in addition that $r'$ lies on a unique component \cite[Definition 3.13]{Tho12}.
\item \emph{Tensor products.} The Rankin--Selberg product of two cuspidal representations of $\GL_2$ over a number field is an isobaric automorphic representation of $\GL_4$, functorial at all places \cite[Theorem M]{Ram00}. For two cuspidal factors without a dihedral factor, this product is cuspidal precisely when the factors are inequivalent up to character twist.
\item \emph{Factors.} \cite[Theorem 1.2]{Die26}: if $\rho_1,\rho_2$ are Lie-irreducible, potentially modular two-dimensional representations of $G_{F'}$, $F'$ totally real, each over its own finite Galois totally real extension, and $\rho_1\otimes\rho_2$ corresponds to a regular L-algebraic cuspidal automorphic representation of $\GL_4(\A_{F'})$, then $\rho_1$ and $\rho_2$ are modular.
\item \emph{Adequacy.} For $p\ge 2(n+1)$, a finite subgroup of $\GL_n(\overline\F_p)$ acting absolutely irreducibly is adequate \cite[\S2 and Appendix]{Tho12}; we use this for $n=2$. From \cite[\S2.2]{ADP}, for finite groups: a subgroup $G\subset\GL(V)$ is \emph{adequate$^+$} if $\mathrm{Ext}^1_G(V,V)=0$, $H^1(G,k)=H^2(G,k)=0$ and $\mathrm{End}(V)$ is spanned by the semisimple elements of $G$; this is equivalent to: $G$ is adequate and $H^2(G,k)=0$ \cite[Lemma 2.6]{ADP}. The image in $\GL(V_1\otimes V_2)$ of a direct product of adequate$^+$ subgroups is adequate$^+$ \cite[Proposition 2.9]{ADP}; a group containing an adequate$^+$ normal subgroup of index prime to $p$ is adequate$^+$ \cite[Lemma 2.10]{ADP}; and an absolutely irreducible subgroup of order prime to $p$ is adequate$^+$ \cite[Lemma 2.11]{ADP}. If $G$ contains a conjugate of $\mathrm{SL}_n(\F_q)$ as a normal subgroup and is contained in its normalizer, the index is prime to $p$ \cite[Proposition 2.12]{ADP}.
\item \emph{Local rings.} The potentially crystalline lifting rings with regular Hodge--Tate weights have formally smooth generic fibre \cite[Theorem 3.3.8]{Kis08}; completed tensor products over $\cO$ of geometrically integral $\cO$-flat complete local rings are integral \cite[Lemma 3.3]{BLGHT}.
\item \emph{Soluble extensions and descent.} \cite[Lemma 4.1.2]{CHT}, also \cite[Lemma A.2.1]{BLGGT}: given a finite Galois extension $D/F$, a finite set of places and finite Galois extensions of the corresponding completions, there is a finite soluble Galois extension of $F$, linearly disjoint from $D$, with those completions. We use soluble automorphic base change and cyclic descent, including the character twist description of the fibres of cyclic base change \cite[\S2.3]{Ram00}, and the polarized formulation of soluble descent \cite[Lemma 2.2.2]{BLGGT}. For the passage from a CM field to its totally real subfield we use cyclic descent and the associated Galois representations \cite[Theorem 2.1.1]{BLGGT}, with the character ambiguity resolved as in the proof of Proposition~\ref{prop:swap} below.
\item \emph{Semistable lifting.} \cite[Theorem 3.6]{ADP} extends (I4) to potentially semistable representations. Let $E_0$ be an imaginary CM field with $\zeta_p\notin E_0$, and let $(r,\nu)$ be an algebraic polarized $p$-adic representation, with $r$ irreducible and $\bar r(G_{E_0(\zeta_p)})$ adequate. Assume that the residual polarized pair $(\bar r,\bar\nu)$ is automorphic, arising from a regular algebraic, cuspidal, polarized representation $(\pi,\chi)$, and that $r|_{G_{E_{0,w}}}\sim r_p(\pi)|_{G_{E_{0,w}}}$ at every finite place $w$. Here $\sim$ is the relation of \cite[Definition 3.27]{ADP}: away from $p$ the representations lie on a common component of the lifting ring, and above $p$ they become semistable over a finite extension, with the same Hodge--Tate weights, on a common component of the corresponding semistable ring. Then $(r,\nu)$ is automorphic. The multiplier of the residual seed is required to be $\bar\nu$. This is the raw polarized multiplier of \cite{ADP}, not the normalized character $\mu$ in (I4). No condition of potential unramifiedness above $p$ is imposed on $\pi$. This is used only in Section~\ref{sec:ss}.
\end{enumerate}

\subsection{The results of \cite{MF} when $p$ ramifies in $F$}\label{ss:ramified}
The results of \cite{MF} are stated for $p$ unramified in $F$. We use the following ones: the twisting Lemma 2.1; the crystalline rings of \S2.2 and the potentially Barsotti--Tate rings and the cycles $\cC_\sigma$ of non-Steinberg Serre weights of \S2.3, which are introduced there for $K/\Q_p$ unramified; the Steinberg ring of \S2.4 at places not above $p$; Definition 3.1 of a global setting, \emph{with the condition that $p$ be unramified in $F$ deleted}; Lemmas 3.2 and 3.3, Proposition 3.4, Lemmas 3.6 and 3.7; the reduction of \S4.1 and Lemma 4.1; and their analogues for general Hodge and inertial types in \S6. From now on a global setting is understood in this extended sense. Proposition 2.2 of \cite{MF} is not used: only the definition of the cycles $\cC_\sigma$ by \cite[(2.1)]{MF} enters, see below. The proofs of the listed results do not use the assumption, for the following reasons.
\begin{itemize}
\item \emph{Local rings.} The crystalline, potentially crystalline and potentially semistable lifting rings and the regularity of the generic fibres of the first two \cite{Kis08}, the twisting by square roots of unramified characters \cite[Lemma 4.3.1]{EG}, and the Steinberg ring at places prime to $p$ do not involve any condition on the local field above $p$. In \cite[\S2.2]{MF} the dimension of the generic fibre is $d=3+[K:\Q_p]$ for any $K$.
\item \emph{Determinants and the twist $\nu$.} If $F_v$ is ramified, a crystalline character of $G_{F_v}$ with Hodge--Tate weight $k-1$ at every embedding is $\eps^{1-k}$ times an unramified character, since two crystalline characters with the same Hodge--Tate weights differ by an unramified one. Hence the character $\psi$ of \cite[\S4.1]{MF} is unramified above $p$, as there. In the proof of \cite[Lemma 3.2]{MF}, an element $z\in\cO_{F_v}^\times$ acts on $\bigotimes_{\iota\colon F_v\to E}\mathrm{Sym}^{k-2}$, the product being over all $[F_v:\Q_p]$ embeddings, by $\prod_\iota\iota(z)^{k-2}=N_{F_v/\Q_p}(z)^{k-2}$, and $\eps(\operatorname{Art}(z))=N_{F_v/\Q_p}(z)^{-1}$; so the computation is unchanged. The character $\psi_2|_{I_{F_v}}$ has order prime to $p$, so it factors through $k_v^\times$ and is the determinant of tame types, as in \cite[\S2.3]{MF}.
\item \emph{Cycles of Serre weights.} The integers $n_\tau(\sigma)$ of \cite[(2.1)]{MF} come from \cite[Lemma 5.1]{CEGS}, which is a statement about representations of $\GL_2(k_v)$. We fix one choice of these integers for each non-Steinberg weight, and use it at all places above $v$ over every extension in which $v$ splits completely. With such a fixed choice, $\cC_\sigma$ is defined by the local rings at $v$ alone, and the proof of \cite[Lemma 3.7(3)]{MF} uses only the patching identity \cite[(3.3)]{MF} and the Grothendieck group identity \cite[(2.1)]{MF}; it does not use the effectivity, the independence of the choice, or the identities of \cite[Proposition 2.2]{MF}. For arbitrary finite $K/\Q_p$, \cite[Theorem 5.2]{CEGS} proves the corresponding assertions without fixing the determinant. The unramified character parameter of \cite[Lemma 4.3.1]{EG} gives the fixed determinant formulation as in \cite[proof of Proposition 2.2]{MF}.
\item \emph{Coefficient systems.} The types $\sigma(\tau_v)$ and $\sigma(\lambda_v,\tau_v)$ are those of \cite[\S2.1]{GK}, defined for any finite extension of $\Q_p$; their semisimplified reductions factor through $\GL_2(k_v)$, since the kernel of reduction is a normal pro-$p$ group and acts trivially on each irreducible constituent. The set $W(\bar r)$ is the set $W^{\mathrm{BT}}(\bar r)$ of \cite{GK}; its explicit description \cite{GLS} is not used.
\item \emph{Patching.} The construction of \cite[\S3]{MF} follows \cite[\S2]{Kis09} and \cite[\S5]{GK}, which allow an arbitrary totally real field and $p>2$; in particular \cite[Lemmas 5.2.5 and 5.3.6, Corollary 5.4.5]{GK} hold for any totally real field, under the standing assumptions of \cite[\S\S5.2.3--5.2.4]{GK}, with arbitrary inertial types in the potentially Barsotti--Tate lifting theorem. In the proof of \cite[Proposition 3.4]{MF} the assumption enters only through the condition of \cite[(2.2)(3)]{Kis09} at $p=5$, which is irrelevant since $p\ge7$; the dimension count, the regularity of the generic fibres and the rank computation there are valid for any local field above $p$.
\par In particular, \cite[\S5.2.3]{GK} uses a residual modular witness of parallel weight two. In the soluble preparation of \S\ref{ss:red} we may take the totally real field to have even degree over $\Q$ by the final quadratic base change of \cite[\S4.1]{MF}, splitting the places above $p$ and in $\Sigma$ and preserving the residual image. Over any such field, transfer a regular Hilbert modular witness to the definite quaternion algebra split at all finite places. Take an integral eigenform whose reduction is nonzero, and shrink its level at $p$ until the reduced coefficient representation is trivial there. In the resulting space of functions on the quaternionic double cosets, a nonzero scalar coordinate of its reduction has the same Hecke eigenvalues away from the level and $p$. Scalar functions at this level form a finite free integral module and reduce onto the scalar functions in characteristic $p$. The Deligne--Serre lifting lemma thus gives a characteristic zero scalar eigensystem with this reduction. A character of the reduced norm would have reducible residual representation. Inverse Jacquet--Langlands gives the required cuspidal eigenform of parallel weight two, with the level at $p$ unrestricted. The degree remains even over $L$, where the same argument applies to the modular witness provided by $\rho|_{G_L}$.
\item \emph{Reduction and the field $L$.} The soluble base change of \cite[\S4.1]{MF} and \cite[Lemma 4.1]{MF} rest on \cite[Lemma 4.1.2]{CHT}, \cite[Lemma 2.2.2]{BLGGT} and \cite[Theorem 2.1.16]{CEG}, none of which assumes $p$ unramified. In the last paragraph of the proof of \cite[Lemma 4.1]{MF}, the assertion that $p$ is unramified in $L$ is to be replaced by the fact that the places above $p$ split completely in $L$, so that the local fields above $p$ are those of $F$; this is all that is used.
\end{itemize}
These arguments extend the comparison needed below to ramified $p$. The digit computations in \cite[\S5]{MF} concern the unramified case.

\section{The auxiliary form}\label{sec:aux}
\begin{defn}\label{def:aux}
Let $\Rb$ be as in Theorem~\ref{thm:main} and let $w$ be an integer and $\mathcal K$ a finite set of integers $\ge2$, all of the parity of $w+1$. An \emph{auxiliary form} for $(\Rb,w,\mathcal K)$ is a Hilbert cuspidal eigenform $g$ over $F$, of parallel weight $k_S\ge2$ with $k_S\notin\mathcal K$ and $k_S\equiv w+1\bmod 2$, of level prime to $p$, whose central character has finite part $\eta$ of finite order, such that, writing $\gt=\rho_g\otimes\eps^{(k_S-1-w)/2}$, whose two Hodge--Tate weights sum to $w$ at every embedding, and $\Sb$ for the reduction of $\gt$:
\begin{enumerate}[label=\textup{(\roman*)}]
\item the projective image of $\Sb$ is not dihedral;
\item the representation $\Rb\otimes\Sb$ of $G_{F(\zeta_p)}$ is absolutely irreducible;
\item the image of $G_{F(\zeta_p)}$ under $\Rb\otimes\Sb$ is adequate.
\end{enumerate}
\end{defn}
For Theorem~\ref{thm:main} we take $w=k-1$ and $\mathcal K=\{k\}$. The congruence $k_S\equiv w+1\bmod2$ makes the required power of $\eps$ integral. The finite character $\eta$ has sign $(-1)^{k_S}$ at every real place; no bound on its order is needed.

\begin{lem}\label{lem:aux}
Assume that $p\ge7$ and that $\Rb|_{G_{F(\zeta_p)}}$ is absolutely irreducible with non-dihedral projective image. Then for every $(w,\mathcal K)$ as above there is an auxiliary form for $(\Rb,w,\mathcal K)$.
\end{lem}
\begin{proof}
Fix $k_S\ge2$ with $k_S\equiv w+1\bmod2$ and $k_S\notin\mathcal K$. Let $\eta_0=1$ if $k_S$ is even. If $k_S$ is odd, let $\eta_0$ be the quadratic character of a CM quadratic extension of $F$ unramified above $p$.

\emph{Existence.} Let $S$ contain the places above $2$ and $p$ and the places where $\eta_0$ is ramified. Let $H/F$ be the compositum of all quadratic extensions unramified at the finite places outside $S$. This is a finite extension. Choose a prime
\[
 t>\max\{5,p,|\Rb(G_F)|\}
\]
which is unramified in $F$ and does not divide any place in $S$. The extension $H/F$ is unramified above $t$, whereas $F(\zeta_t)/F$ is totally ramified there, so $H\cap F(\zeta_t)=F$. By Chebotarev there are primes $\mathfrak q\notin S$, of arbitrarily large norm, which split completely in $H$ and satisfy $N\mathfrak q\equiv-1\bmod t$.

At $\mathfrak q$ choose the supercuspidal type attached to a tame character of the multiplicative group of the unramified quadratic extension of $F_{\mathfrak q}$, trivial on a uniformizer and of order $t$ on the units. The character of the quadratic residue field has niveau two. The associated cuspidal type has dimension $N\mathfrak q-1$ and trivial central character on units of $F_{\mathfrak q}$. At the primes ramified in $\eta_0$ prescribe the principal series type of $1\oplus\eta_0$, and prescribe the trivial type at all other finite places. At infinity prescribe the discrete series of parallel weight $k_S$. These types satisfy the central character compatibility of \cite[\S1.2]{We09}, since $\eta_0$ has the required sign. Their number of nontrivial finite components is bounded, whereas their dimension product tends to infinity with $N\mathfrak q$. By \cite[Theorem 1.1]{We09} they occur in a Hilbert cuspidal eigenform $g$, of level prime to $p$.

The finite part of the central character of $g$ is $\eta=\eta_0\xi$, where $\xi$ is an everywhere unramified character of the ideal class group. In particular $\eta$ has finite order and is unramified above $p$. We take this actual character as the $\eta$ in Definition~\ref{def:aux}.

At $\mathfrak q$ the two tame inertia characters are $\alpha$ and $\alpha^{-1}$ of order $t$. Since $t\ne p$, they remain distinct of order $t$ after reduction. Frobenius interchanges their lines because $N\mathfrak q\equiv -1\pmod t$. This proves local absolute irreducibility, and the quotient $\bar\alpha^2$ still has order $t$ in projective inertia. Thus its projective image is neither reducible nor $A_4$, $S_4$ or $A_5$. If it were dihedral, let $K/F$ be its quadratic inducing field. This field is unramified at $\mathfrak q$, since the inertia image there has odd order, and is unramified at every finite place outside $S\cup\{\mathfrak q\}$. Hence $K\subset H$. But $\mathfrak q$ splits in $H$, so it splits in $K$, contradicting local irreducibility at $\mathfrak q$. This is the argument of \cite[Lemma 4.2]{DZ}, with the fixed residual prime included in $S$. The use of tamely dihedral primes to force a large residual image goes back to \cite{DW}, for classical modular forms; \cite{DZ} extends it to Hilbert modular forms.

By Dickson's classification the projective image $P_S$ is $\PSL_2(\F')$ or $\PGL_2(\F')$. Its normal subgroup $Q_0=\PSL_2(\F')$ contains the element of order $t$, since $[P_S:Q_0]\le2$ and $t$ is odd. Therefore $|Q_0|\ge t>|\Rb(G_F)|$. Twisting by a power of $\eps$ does not change this projective image. This proves (i). Let $P_R$ be the projective image of $\Rb$.

\emph{Irreducibility.} Let $P\subset P_R\times P_S$ be the image of $G_{F(\zeta_p)}$. Its projection to $P_S$ still contains $Q_0$, which is perfect. The kernel of the projection of $P$ to the first factor is $1\times N$ with $N$ normal in the second projection, and the quotient of the second projection by $N$ is a quotient of a subgroup of $P_R$. As $Q_0$ is simple, $N\cap Q_0$ is $Q_0$ or trivial, and the second case would embed $Q_0$ in a group of order at most $|P_R|$. So $P\supseteq1\times Q_0$. Hence the image of $G_{F(\zeta_p)}$ under $(\Rb,\Sb)$ contains a normal subgroup $H_0$ which acts on $\Rb$ by scalars and absolutely irreducibly on $\Sb$. A subspace of $\Rb\otimes\Sb$ stable under $H_0$ is of the form $W\otimes\Sb$, and it is stable under the whole image only if $W$ is stable under $G_{F(\zeta_p)}$. As $\Rb|_{G_{F(\zeta_p)}}$ is absolutely irreducible, so is $\Rb\otimes\Sb$. This gives (ii).

\emph{Adequacy.} We use (I7). Write $F_0=F(\zeta_p)$, $G_R=\Rb(G_{F_0})$ and $G_S=\Sb(G_{F_0})$. We first note that a finite absolutely irreducible subgroup $G\subset\GL_2(\overline\F_p)$ is adequate$^+$. It is adequate because $p\ge7\ge2(2+1)$. And $H^2(G,k)=0$: by Dickson's classification either $G$ has order prime to $p$, or its projective image is $\PSL_2(\F_q)$ or $\PGL_2(\F_q)$ with $q$ a power of $p$; in the second case $G$ contains a conjugate of $\mathrm{SL}_2(\F_q)$ as a normal subgroup of index prime to $p$, the restriction $H^2(G,k)\to H^2(\mathrm{SL}_2(\F_q),k)$ is injective, and the target vanishes by the universal coefficient theorem, since $\mathrm{SL}_2(\F_q)$ is perfect with trivial Schur multiplier for $q\ge7$, $q\neq9$.

The group $G_S$ has projective image $\PSL_2(\F')$ or $\PGL_2(\F')$, so it is contained in the normalizer of $N_S=\mathrm{SL}_2(\F')$ and contains it: the derived subgroup of the inverse image of $\PSL_2(\F')$ is $\mathrm{SL}_2(\F')$. The quotient $G_S/N_S$ is abelian, because it embeds in the quotient of the scalar multiples of $\GL_2(\F')$ by $\mathrm{SL}_2(\F')$, and of order prime to $p$. Let $K_R$, $K_S$ be the fields cut out by $\Rb$, $\Sb$ over $F_0$ and $K=K_R\cap K_S$. The group $\mathrm{Gal}(K/F_0)$ is a common quotient of $G_R$ and $G_S$. The image of $N_S$ in it is a quotient of $\mathrm{SL}_2(\F')$ of order at most $|G_R|<|Q_0|$, hence trivial. So $\mathrm{Gal}(K/F_0)$ is a quotient of $G_S/N_S$: it is abelian of order prime to $p$. The fields $K_R$ and $K_S$ are linearly disjoint over $K$, so the image of $G_K$ under $(\Rb,\Sb)$ is $\Rb(G_K)\times\Sb(G_K)$. The group $\Sb(G_K)$ contains $N_S$, so it is absolutely irreducible. The group $\Rb(G_K)$ is normal in $G_R$ with abelian quotient, so it contains the derived subgroup of $G_R$. The projective image of $G_R$ is $A_4$, $S_4$, $A_5$, $\PSL_2(\F_q)$ or $\PGL_2(\F_q)$, by the non-dihedral hypothesis; its derived subgroup is the Klein four-group, $A_4$, $A_5$, $\PSL_2(\F_q)$ or $\PSL_2(\F_q)$, and a subgroup of $\GL_2$ whose projective image contains the Klein four-group is absolutely irreducible. So $\Rb(G_K)$ is absolutely irreducible. By the first paragraph $\Rb(G_K)$ and $\Sb(G_K)$ are adequate$^+$. By \cite[Proposition 2.9]{ADP} the image of $G_K$ under $\Rb\otimes\Sb$ is adequate$^+$, and by \cite[Lemma 2.10]{ADP} so is the image of $G_{F_0}$, since $[K:F_0]$ is prime to $p$. As $p$ does not divide $4$, it is adequate in the sense of \cite{Tho12}. This is the argument of \cite[Remark 2.13]{ADP}, and gives (iii).
\end{proof}

\begin{rem}[The fixed residual prime]\label{rem:NT}
The construction above requires a large image only in the chosen characteristic $p$, so one auxiliary place suffices. To control all residual characteristics at once, \cite[Proposition 6.1]{DZ} uses two tamely dihedral places with distinct residue characteristics and distinct auxiliary prime orders. An analogous result for forms of parallel weight $2$ is \cite[Proposition 5.4]{NT}. Including the places above the fixed $p$ in $S$ avoids any restriction on the weight in the preceding argument.
\end{rem}

\section{Fields and local components}\label{sec:setup}
In this section and the next two we prove Theorem~\ref{thm:main}; so $\rho$, $\Rb$, $k$ are as there.

\subsection{Reduction}\label{ss:red} By \cite[\S4.1]{MF} we may replace $F$ by a soluble totally real extension in which the places above $p$ split completely, and assume that we are in a global setting $(F,\Sigma,\bar r,\psi)$ of weight $k$ with $\bar r=\Rb$, in which $\rho$ has the properties $(\dagger)$. Hypothesis (3) is preserved: the extension is linearly disjoint from the field cut out by $\Rb$ and $\zeta_p$. If $\rho$ is modular over the new field, it is modular over the original one by soluble descent.

Let $\cA$ be the set of automorphic components of the setting. For $C=(C_v)_{v\mid p}\in\cA$, fix an eigenform contributing to $M_k$ and write $f_C$ for its untwisted Hilbert Galois representation. If $\nu_{\mathrm{MF}}$ is the character denoted $\nu$ in \cite[\S3.1]{MF}, then $(f_C\otimes\nu_{\mathrm{MF}})|_{G_{F_v}}$ lies on $C_v$. In the reduction of \cite[\S4.1]{MF}, $\psi=\det\rho\,\eps^{k-1}$, so $\det f_C=\theta_k=\det\rho$; the representation $f_C$ has the properties $(\dagger)$. By Lemma~\ref{lem:aux} fix an auxiliary form $g$ for $(\Rb,k-1,\{k\})$, with $\gt$, $\Sb$ and $k_S$ as in Definition~\ref{def:aux}.

\subsection{The field $L$} Apply \cite[Theorem 2.1.16]{CEG}, as in the proof of \cite[Lemma 4.1]{MF}, to the compatible systems of $\rho$, of the $f_C$ for $C\in\cA$, and of $g$, taking the field to avoid to contain the fields cut out by $\Rb$, $\Sb$ and $\zeta_p$. The argument of \cite[Lemma 4.1]{MF} does not use that the representations lift the same residual representation: each system is odd, regular, pure and polarized by its determinant, and becomes modular over a totally real field. We get a totally real Galois extension $L/F$ in which the places of $F$ above $p$ and in $\Sigma$ split completely, linearly disjoint over $F$ from the field cut out by $\Rb$, $\Sb$ and $\zeta_p$, and over which $\rho$, all the $f_C$ and $g$ are modular. The construction does not ensure that $L/F$ is soluble. Thus the modularity of the $f_C$ and of $g$ over $L$ is not a consequence of their modularity over $F$. The base change to $L$ of the global setting is a global setting of weight $k$ over $L$, with set of automorphic components $\cA_L$ \cite[Lemma 4.1]{MF}. Let $d=[L:F]$. We fix a partition of the places of $L$ above $p$ into $d$ \emph{sheets}, each sheet containing exactly one place above each $v\mid p$ of $F$.

\subsection{The field $F'$} Let $D$ be the compositum of $L$ with the fields cut out by $\Rb$, by $\Sb$ and by $\zeta_p$; it is Galois over $F$. By (I9) there is a finite soluble Galois totally real extension $F'/F$, linearly disjoint from $D$, such that:
\begin{itemize}
\item for every place $w$ of $F'$ above $p$, the restrictions of $\Rb$ and $\Sb$ to $G_{F'_w}$ are trivial, and $\det\rho$ and $\det\gt$ agree on $G_{F'_w}$;
\item over $F'$ the form $g$ is unramified outside $p$, and $\Sb|_{G_{F'_w}}$ is trivial for $w$ above $\Sigma$.
\end{itemize}
For the first point, $\det\rho/\det\gt$ has finite order, by global class field theory. Indeed their quotient has Hodge--Tate weight zero and finite inertia at every place above $p$, and has finite inertia at the other finite places. It therefore factors through a finite ray class group. For the second, the local components of $g$ at its primes of ramification are induced from tamely ramified characters at $\mathfrak q$, and ramified principal series at the primes dividing the conductor of $\eta$; they become unramified over a finite extension of the local field. The real places are required to split. The extension $F'/F$ is ramified above $p$ in general.

Put $L'=LF'$. By soluble base change $\rho$, the $f_C$ and $\gt$ are modular over $L'$, and every place of $L'$ above $p$ lies over a place of $L$ in a well-defined sheet. The field cut out by $\Rb$, $\Sb$ and $\zeta_p$ is linearly disjoint from $L'$ over $F$. For each $v\mid p$ fix a completion $K'_v$ of $F'$ at $v$. For each $u\mid v$ in $L'$ choose a prolongation and an $F_v$-isomorphism $L'_u\simeq K'_v$ compatible with the chosen copy of $F_v=L_{u|L}$. Such identifications exist since $v$ splits completely in $L$ and $F'/F$ is Galois. We use them in restricting representations at $u$. We write $\chi_R=\det\rho|_{G_{L'}}$ and $\chi_S=\det\gt|_{G_{L'}}$; they agree on the decomposition groups above $p$.

\subsection{Potentially connected representations}\label{ss:potconn}
\begin{defn}\label{def:potconn}
Let $K$ be a finite extension of $\Q_p$ and let $r_1,r_2\colon G_K\to\GL_n(\cO)$ be continuous representations. We say that $r_1$ and $r_2$ are \emph{potentially connected} if there is a finite extension $K'/K$ such that:
\begin{enumerate}[label=\textup{(\roman*)}]
\item the reductions of $r_1|_{G_{K'}}$ and $r_2|_{G_{K'}}$ are trivial;
\item $r_1|_{G_{K'}}$ and $r_2|_{G_{K'}}$ are crystalline, with the same Hodge--Tate weights at every embedding, and these are distinct;
\item $r_1|_{G_{K'}}$ and $r_2|_{G_{K'}}$ define points of the same irreducible component of the generic fibre of the crystalline framed deformation ring of the trivial representation of $G_{K'}$ with these Hodge--Tate weights.
\end{enumerate}
Two representations of the Galois group of a number field are \emph{potentially connected above $p$} if their restrictions to the decomposition group at every place above $p$ are potentially connected.
\end{defn}

\begin{rem}\label{rem:potconn}
(a) The residual representations of $r_1$ and $r_2$ over $K$ need not be isomorphic: only their restrictions to $G_{K'}$ are compared. This is the point of the definition.

(b) If the conditions hold for $K'$ they hold for every finite extension $K''$ of $K'$: the restrictions to $G_{K''}$ of the points of a component over $K'$ form an irreducible set, which lies on one component over $K''$. Hence the relation is an equivalence relation on the potentially crystalline representations with given regular Hodge--Tate weights.

(c) By the second part of the proof of Lemma~\ref{lem:tensor}, which applies to the crystalline deformation ring of a trivial residual representation in any dimension, a component of the ring in (iii) is stable under conjugation by $\GL_n(\cO)$. So the relation does not change if $r_1$ or $r_2$ is replaced by a conjugate under $\GL_n(\cO)$.

(d) If $r_1$ and $r_2$ have the same reduction and lie on the same component of a potentially crystalline lifting ring over $K$, that is $r_1\sim r_2$ in the sense of \cite{BLGGT} and \cite[Definition 3.7]{Tho12}, they are potentially connected. A representation is potentially diagonalizable in the sense of \cite{BLGGT} when it is potentially connected to a sum of crystalline characters. Conversely, if $r_1$ and $r_2$ are potentially connected, their restrictions to $G_{K'}$ satisfy $r_1|_{G_{K'}}\sim r_2|_{G_{K'}}$, which is the local hypothesis above $p$ of \cite[Theorem 7.1]{Tho12} over a field whose completions are the $K'$.
\end{rem}

\subsection{Local components over $K'_v$}\label{ss:local}
For $v\mid p$ consider the following reduced, $\cO$-flat framed deformation rings of the trivial representation of $G_{K'_v}$, with determinant $\chi_R|_{G_{K'_v}}=\chi_S|_{G_{K'_v}}$ for the first two:
\begin{itemize}
\item $R^{k}_{v}$: crystalline lifts with Hodge--Tate weights $0,k-1$ at every embedding;
\item $R^{g}_{v}$: crystalline lifts with the Hodge--Tate weights of $\gt$ at every embedding;
\item $R^{(4)}_{v}$: crystalline lifts of the trivial four-dimensional representation with the Hodge--Tate weights of $\rho\otimes\gt$ at every embedding, and determinant $(\chi_R\chi_S)^2$.
\end{itemize}
The four Hodge--Tate weights are distinct because $k\neq k_S$ and $k,k_S\ge2$. For fixed Hodge and inertial types the determinant varies by an unramified character. The unique $n$th root, congruent to $1$, of an unramified character lifting the trivial one exists for $n=2,4$ since $p\ge7$. Twisting by it identifies the ring with variable determinant with a formal power series ring in one variable over its fixed determinant quotient. Thus (I8) gives formal smoothness for all three generic fibres, so every point lies on exactly one irreducible component. After enlarging the coefficients, the components are geometrically irreducible.

First transport a component $C_v$ of the crystalline ring in \cite[\S2.2]{MF} by the twist $\nu_{\mathrm{MF}}^{-1}$, and then restrict its points to $G_{K'_v}$. Their determinant is $\chi_R|_{G_{K'_v}}$. We use embeddings of the completions over $F_v$ to make these restrictions. Changing such an embedding conjugates the restricted representation of $f_C$ by an integral matrix, and likewise for $\gt$. The components containing these points are unchanged by Remark~\ref{rem:potconn}(c). The image of the component under restriction is irreducible; let $C'_v$ be the unique component of $R^k_v$ containing it. In particular $f_C|_{G_{K'_v}}$ lies on $C'_v$. Let $E_v$ be the component of $R^g_v$ containing $\gt|_{G_{K'_v}}$. A \emph{letter} is either a family $C'=(C'_v)_v$ with $C\in\cA$, of \emph{weight} $k$, or the family $E=(E_v)_v$, of \emph{weight} $k_S$. For a point of $R^k_v$ (respectively $R^g_v$), membership in $C'_v$ (respectively $E_v$) is equivalent to lying on the same crystalline component over $K'_v$ as the corresponding local representation of $f_C$ (respectively $\gt$). Indeed the fixed determinant condition removes only the unramified parameter described above.

\begin{lem}\label{lem:tensor}
Let $X_v$ be a component of $R^k_v$ and $Y_v$ a component of $R^g_v$. There is a unique irreducible component $T(X_v,Y_v)$ of $R^{(4)}_v[1/p]$ containing $x\otimes y$ for all points $x$ of $X_v$ and $y$ of $Y_v$. It also contains $y\otimes x$ for all such points.
\end{lem}
\begin{proof}
The Kronecker product of lifts gives a morphism from the spectrum of the completed tensor product of the two quotient rings to $\operatorname{Spec}R^{(4)}_v$: a tensor product of crystalline representations is crystalline with the sums of the weights. The map first takes values in the unrestricted framed lifting space. The pairs of characteristic zero points give a dense set in the spectrum of the completed tensor product by \cite[Lemma 3.3]{BLGHT}, so their Hodge type and determinant make the map factor through the displayed quotient. The source is integral by (I8), so its image is irreducible, and it lies on a single component because the components of the generic fibre are disjoint.

Let $P\in\GL_4(\cO)$ be the permutation matrix with $P(x\otimes y)P^{-1}=y\otimes x$. We prove that every component of a framed crystalline ring $R$ of a trivial residual representation in dimension $n$ is stable under conjugation by $\GL_n(\cO)$, also when its determinant is fixed. Write $R^\square$ for the corresponding unrestricted framed ring and $\mathfrak m$ for the maximal ideal of $R$.

\emph{Elementary conjugations.} Fix a minimal prime $\mathfrak p$ of $R$, put $A=R/\mathfrak p$, and let $B=A\langle s\rangle$ be the $\mathfrak m$-adic completion of $A[s]$. It consists of the series $\sum_{i\ge0}c_i s^i$ with $c_i\to0$, and is a domain because it embeds in $A[[s]]$. Conjugate the universal representation over $A$ by $e_{ij}(s)$. This family is classified by a homomorphism $\phi:R^\square\to B$. Indeed, modulo $\mathfrak m^aB$ its matrix entries minus the identity have bounded polynomial degree, take finitely many values, and have coefficients in the nilpotent ideal $\mathfrak m/\mathfrak m^a$. They generate a finite local $\cO$-subalgebra of $B/\mathfrak m^aB$ with the prescribed residue field. The compatible classifying maps to these algebras give $\phi$.

For every point $x:A\to\Qbar_p$ and every $b\in\cO$, specialization at $(x,s=b)$ is conjugate to $x$ and is crystalline of the prescribed Hodge type. Thus if $f\in\ker(R^\square\to R)$ and $\phi(f)=\sum_i c_i s^i$, the series $\sum_i x(c_i)s^i$ vanishes for every $s\in\cO$. Strassmann's theorem makes it zero, so $x(c_i)=0$ for all $x$ and $i$. The characteristic zero points of the $\cO$-flat domain $A$ are dense \cite[Lemma 3.3(1)]{BLGHT}, whence $c_i=0$ and $\phi$ factors through $R$. Its kernel is prime and is contained in $\mathfrak p$, by specialization at $s=0$, so it equals $\mathfrak p$. Specialization at $s=b$ and then at $s=-b$ proves that conjugation by $e_{ij}(b)$ preserves $\mathfrak p$.

\emph{General conjugations.} For $h\in\GL_n(\cO)$ choose a finite coefficient extension $\cO'/\cO$ and a unit $u\in\cO'$ with $u^n=\det h$. Conjugation by $h$ equals conjugation by $u^{-1}h\in\mathrm{SL}_n(\cO')$, which is a product of elementary matrices. Apply the preceding argument to $R\otimes_\cO\cO'$. Every minimal prime of this ring is stable. Choose a minimal prime above $\mathfrak p$; faithful flatness makes its contraction $\mathfrak p$, and compatibility with conjugation gives $h(\mathfrak p)=\mathfrak p$. In particular conjugation by $P$ preserves the component containing the products $x\otimes y$, so it also contains all $y\otimes x$.

\end{proof}

\begin{cor}\label{cor:potconn}
Let $x,x'$ be points of the same component of $R^k_v$ and $y,y'$ points of the same component of $R^g_v$. Then the representations $x\otimes y$, $x'\otimes y'$ and $y'\otimes x'$ of $G_{K'_v}$ are pairwise potentially connected.
\end{cor}
\begin{proof}
They are points of the component of Lemma~\ref{lem:tensor}, and the conditions of Definition~\ref{def:potconn} hold with $K'=K'_v$.
\end{proof}

\section{Mixed lifts}\label{sec:mixed}
\begin{defn}
A \emph{pattern} is a map $T$ from the set of sheets to the set of letters. Let $\bar r_0$ be $\Rb|_{G_{L'}}$ or $\Sb|_{G_{L'}}$. A \emph{lift of $\bar r_0$ of pattern $T$} is a lift $A\colon G_{L'}\to\GL_2(\Qbar_p)$ of $\bar r_0$, with determinant $\chi_R$ if $\bar r_0=\Rb|_{G_{L'}}$ and $\chi_S$ if $\bar r_0=\Sb|_{G_{L'}}$, such that, at every place $u\mid p$ of $L'$, lying over $v$ and in the sheet $i$, the restriction $A|_{G_{L'_u}}$ is a point of the component $T(i)_v$; and such that, away from $p$: if $\bar r_0=\Rb|_{G_{L'}}$, $A$ is unramified outside $\Sigma$ and corresponds to a point of the Steinberg ring of \cite[\S2.4]{MF} at the places above $\Sigma$; if $\bar r_0=\Sb|_{G_{L'}}$, $A$ is unramified.
\end{defn}
 Here and below the Steinberg ring for an untwisted lift over a completion $K$ above $\Sigma$ is the ring of \cite[\S2.4]{MF} formed over $K$ with determinant $(\det\rho)|_{G_K}$. Twisting by the unramified character $\nu_{\mathrm{MF}}|_{G_K}$ identifies it with the ring for the determinant $(\det\rho)\nu_{\mathrm{MF}}^2$ used in patching. This twist preserves nonzero unipotent monodromy and the ratio of the Frobenius eigenvalues at unramified points. Restriction to a finite extension carries Steinberg points to Steinberg points and carries their closure into the corresponding closure over the extension. The restrictions to $G_{L'}$ of the $f_C$ and of $\gt$ are lifts of the constant patterns $C'$ and $E$.

\begin{lem}\label{lem:exist}
For every pattern $T$ and each of the two residual representations there is a lift of pattern $T$. Every lift of pattern $T$ is Lie-irreducible and potentially modular, over a finite Galois totally real extension of $L'$.
\end{lem}
\begin{proof}
Existence is (I2) over $L'$: the residual representations are totally odd and irreducible over $L'(\zeta_p)$ by the disjointness in \S\ref{sec:setup}; the determinant is $\chi_R$ or $\chi_S$; the components above $p$ are the $T(i)_v$. In the notation of (I2) we take $E_v=K'_v$: the ring there is the ring of semistable lifts of the trivial representation of $G_{K'_v}$ with the given Hodge--Tate weights and determinant $\chi_R|=\chi_S|$. It is non-zero, and the crystalline ring $R^k_v$ or $R^g_v$ is a quotient of it whose generic fibre is equidimensional of the same dimension \cite[Theorems 3.3.4 and 3.3.8]{Kis08}, so the $T(i)_v$ are irreducible components of its generic fibre. At the primes above $\Sigma$, the Steinberg ring selects a component of the full fixed determinant ring: its generic fibre has dimension $3$ by \cite[Proposition 2.4]{MF}, the dimension of every component of the framed generic fibre with fixed determinant away from $p$ \cite[\S1.3]{BLGGT}. The set of places consists of the places above $p$ and, for $\Rb$, those above $\Sigma$, where the condition is the Steinberg component; $\Rb|_{G_{L'}}$ and $\chi_R$ are unramified outside these places by the reduction of \S\ref{ss:red}, and $\Sb|_{G_{L'}}$ and $\chi_S$ are unramified outside $p$ by the choice of $F'$.

Every lift $A$ is irreducible because its reduction is. If its restriction to an open normal subgroup were reducible, Clifford theory gives either two distinct characters interchanged by the whole group or a scalar restriction. The scalar case is incompatible with distinct Hodge--Tate weights. In the other case $A$ is induced from a quadratic extension, and its semisimplified reduction is induced from the reduction of that character. This contradicts the preserved non-dihedral residual hypothesis. Thus $A$ is Lie-irreducible.

Let $M$ be an imaginary quadratic field in which $p$ splits and which is not contained in the compositum of $L'$ with the fields cut out by $\Rb$, $\Sb$ and $\zeta_p$; then $\zeta_p\notin E_0:=L'M$, and $\Rb$, $\Sb$ and $\Rb\otimes\Sb$ have the same images on $G_{E_0(\zeta_p)}$ as on $G_{F(\zeta_p)}$. The restriction $s=A|_{G_{E_0}}$ is polarized by its determinant, odd, and reasonable: $p\ge7$, $\zeta_p\notin E_0$, and its reduction is irreducible on $G_{E_0(\zeta_p)}$ by the choice of $M$. The places of $E_0$ above $p$ are split over $L'$, because $p$ splits in $M$, and their completions are the $K'_v$; so the local polarized ring in (I3) is the $\GL_2$ lifting ring with variable determinant \cite[Remark 2.1.11]{CEG}. Each fixed determinant component $C'_v$ or $E_v$ has irreducible image in the full potentially semistable space of trivial inertial type. Choose a full component containing that image. It contains both the point of $A$ at the place in question and the modular point defining the letter. The latter witnesses global realizability of this full component by \cite[Definition 2.1.9 and Lemma 2.1.15]{CEG}. For $C'_v$ take the compatible system of $f_C|_{G_{E_0}}$, and for $E_v$ that of $\gt|_{G_{E_0}}$. These are odd, regular, polarized and weakly irreducible, because $f_C$ and $g$ have no complex multiplication; their reductions $\Rb|_{G_{E_0}}$ and $\Sb|_{G_{E_0}}$ are reasonable; at a place above $v$ their local residual representations are trivial, like that of $A$; and their local representations lie on $C'_v$ and $E_v$ by definition. The multiplier plays no role, by the splitting above. By (I3), $s$ is automorphic over a finite Galois CM extension $N/E_0$. A CM field is stable under complex conjugation, which is central in its Galois group, so $N$ is Galois over $L'$ and its maximal totally real subfield $N^+$ is Galois over $L'$. The restriction $A|_{G_N}$ is irreducible by the preceding paragraph. Cyclic automorphic descent along $N/N^+$, with the character ambiguity removed as in Proposition~\ref{prop:swap}, shows that $A|_{G_{N^+}}$ is modular. In rank two the Hecke avatar of $(\det A)\eps$ supplies the essential self-duality character in the C-algebraic normalization. Thus $A$ is potentially modular over the finite Galois totally real extension $N^+/L'$.

\end{proof}

\begin{prop}[Tensor swap]\label{prop:swap}
Let $A,A_0$ be lifts of $\Rb|_{G_{L'}}$ of patterns $T,T_0$, and $B,B_0$ lifts of $\Sb|_{G_{L'}}$ of patterns $U,U_0$. Assume that for every sheet $i$ the letters $T(i)$ and $U(i)$ have different weights, and $\{T(i),U(i)\}=\{T_0(i),U_0(i)\}$. If $A_0$ and $B_0$ are modular, then $A$ and $B$ are modular.
\end{prop}
\begin{proof}
Let $\pi_A,\pi_B$ be the L-algebraic cuspidal representations attached to $A_0,B_0$ in our Satake normalization, and put $\Pi_0=\pi_A\boxtimes\pi_B$. It corresponds to $A_0\otimes B_0$ and is regular L-algebraic, since the four Hodge--Tate weights are distinct. Put $\Omega_0=\Pi_0|\det|^{3/2}$. Then $\Omega_0$ is regular C-algebraic and essentially self-dual, and, once cuspidality is verified below, $r_p(\Omega_0)=A_0\otimes B_0$ by the normalization in \cite[Theorem 2.1.1]{BLGGT}.

 The cuspidality criterion of \cite[Theorem M]{Ram00} applies directly. Neither factor is dihedral by Lemma~\ref{lem:exist}. At each archimedean embedding the differences of the two weights for the factors are distinct. A character twist preserves this difference, so the two factors are not twist equivalent. Hence $\Pi_0$ is cuspidal.

Let $M$ and $E_0=L'M$ be as in the proof of Lemma~\ref{lem:exist}. The two rank two base changes are still cuspidal and non-dihedral by Lie irreducibility and cyclic base change. Their archimedean weight differences remain distinct. The same criterion gives a cuspidal product over $E_0$, which is the base change of $\Pi_0$ by its unramified parameters. Write $r=(A\otimes B)|_{G_{E_0}}$, $r'=(A_0\otimes B_0)|_{G_{E_0}}$, $m=\chi_R\chi_S$ and $\xi=m^{-1}\eps^{-3}$. Choose an algebraic character $\theta$ of $G_{E_0}$, crystalline at $p$, with $\theta\theta^c=\xi|_{G_{E_0}}$, by \cite[Lemma A.2.5(1)]{BLGGT}. Indeed every place above $p$ splits over $L'$. In each conjugate pair prescribe the trivial character at one place and the restriction of $\xi$ at the other. These local characters are crystalline, their product after conjugation is $\xi$ on inertia, and $\xi(c)=-1$ at every real place. Agreement on inertia gives a crystalline character at each place above $p$. Put $r_\theta=r\otimes\theta$ and $r'_\theta=r'\otimes\theta$, and write $\Theta$ for the corresponding Hecke character. We apply (I4) to these twists.
\begin{itemize}
\item \emph{Polarization.}  Both untwisted tensors have a symmetric pairing over $L'$ with multiplier $m$, obtained from the alternating determinant pairing on each factor. Thus $r_\theta^c\simeq r_\theta^\vee\eps^{-3}$ and likewise for $r'_\theta$. At almost all places local compatibility and $\theta\theta^c=m^{-1}\eps^{-3}$ show that the C-algebraic seed $\mathrm{BC}_{E_0/L'}(\Omega_0)\otimes\Theta$ is conjugate self-dual. Strong multiplicity one gives this identity globally. Its normalizing Hecke character in (I4) is $1$, whose archimedean sign is $+1$ as required for $n=4$. To check the character on $G_{L'}$, write $\delta=\delta_{E_0/L'}$. The symmetric pairing gives the odd raw multiplier $m\delta$ before twisting \cite[\S1.1]{BLGGT}. After twisting it is $m\delta(\theta\circ\operatorname{Ver}_{E_0/L'})$. This agrees with $\eps^{-3}$ on $G_{E_0}$ by the norm identity and at complex conjugation $c$ because transfer sends $c$ to $c^2=1$. In bases identifying the residual factors, their determinant pairings and hence the residual tensor pairings agree. This verifies the congruence of polarized pairs.
\item \emph{Residual image.}  The untwisted representations $r$ and $r'$ reduce to $(\Rb\otimes\Sb)|_{G_{E_0}}$, which is absolutely irreducible on $G_{E_0(\zeta_p)}$ with adequate image, by Definition~\ref{def:aux}(ii), (iii) and disjointness. The common scalar twist also has adequate image. In the joint image of $(\bar r,\bar\theta)$, the maps to the images of $\bar r$ and $\bar r\otimes\bar\theta$ have central kernels of order prime to $p$ acting trivially on the adjoint module. Inflation thus identifies the cohomology groups in the adequacy conditions. Multiplication by a scalar of order prime to $p$ preserves semisimplicity and changes the spanning matrices only by nonzero scalar factors.
\item \emph{Above $p$.} The base change of $\Pi_0$ is unramified above $p$, because $A_0$ and $B_0$ are crystalline over the $K'_v$. At a place over $v$ in the sheet $i$, both $r$ and $r'$ are crystalline over $K'_v$ with the same weights. If $T(i)=T_0(i)$, the two local representations are both of the form $x\otimes y$, with $x$ in the letter of weight $k$ and $y$ in $E_v$, or both of the form $y\otimes x$. If $T(i)=U_0(i)$, one is of the form $x\otimes y$ and the other of the form $y\otimes x$. In all cases Corollary~\ref{cor:potconn} shows that they are potentially connected, through $K'_v$ itself: they lie on a common component of $R^{(4)}_v[1/p]$. So $A\otimes B$ and $A_0\otimes B_0$ are potentially connected above $p$, and over $E_0$, whose completions above $p$ are the $K'_v$, this is $r'\sim r$. Twisting by the fixed crystalline character $\theta$ identifies the lifting rings after shifting each labelled Hodge weight by the weight of $\theta$. The two twisted points therefore lie on a common component, and the twisted seed is still unramified above $p$.
\item \emph{Away from $p$.} At a place not above $\Sigma$, all four of $A,B,A_0,B_0$ are unramified, so $r$ and the automorphic representation are unramified. At a place $w$ above $\Sigma$, the residual representations are trivial, $B$ and $B_0$ are unramified, and $A$, $A_0$ are points of the Steinberg ring, which is a domain; the ring of unramified lifts is formally smooth. As in Lemma~\ref{lem:tensor}, the tensor products of the points of these two rings lie on one irreducible component of the unrestricted lifting ring at $w$, so $r'\sim r$. The local component at $w$ of the cuspidal representation attached to $r'$ is generic, so $r'$ is a smooth point by \cite[Lemma 1.3.2(1)]{BLGGT}; hence it lies on a unique component and $r'\leadsto r$. Twisting by $\theta$ identifies the unrestricted lifting rings at each such place. If $\theta$ introduces ramification at a place where both tensors were unramified, the two original points lie on the common unramified component \cite[\S1.3]{BLGGT}. After twisting, the seed is still generic and its point is smooth. This gives the required strong connection also at these places.
\end{itemize}
 The lifting theorem gives a regular C-algebraic cuspidal representation over $E_0$ attached to $r_\theta$. Twist it by $\Theta^{-1}$ and denote the result by $\Omega_E$. Then $r_p(\Omega_E)\simeq(A\otimes B)|_{G_{E_0}}$. Since $A\otimes B$ extends to $G_{L'}$, Chebotarev and strong multiplicity one give $\Omega_E^c\simeq\Omega_E$. Cyclic descent gives a cuspidal representation $\Omega$ over $L'$, whose archimedean parameters make it regular C-algebraic. Let $\delta=\delta_{E_0/L'}$ and let $\chi$ be the algebraic Hecke character with Galois avatar $\chi_R\chi_S\eps^3$. Local compatibility gives $\Omega_E^\vee\otimes\chi_{E_0}\simeq\Omega_E$. The fibre of cyclic base change then gives
\[
\Omega^\vee\otimes\chi\simeq\Omega\otimes\delta^e,\qquad e\in\{0,1\}.
\]
Hence $\Omega$ is essentially self-dual and has an associated Galois representation \cite[Theorem 2.1.1]{BLGGT}. Its restriction to $G_{E_0}$ is the irreducible representation $(A\otimes B)|_{G_{E_0}}$, so Clifford theory gives $r_p(\Omega)\simeq(A\otimes B)\otimes\delta^a$ for some $a\in\{0,1\}$. Replacing $\Omega$ by $\Omega\otimes\delta^a$ makes $r_p(\Omega)\simeq A\otimes B$. Therefore $\Omega|\det|^{-3/2}$ is the regular L-algebraic cuspidal representation required by (I6). Lemma~\ref{lem:exist} and (I6) imply that $A$ and $B$ are modular.
\end{proof}

\begin{prop}\label{prop:mixed}
Assume that $\cA$ is not empty. For every pattern $T$, every lift of $\Rb|_{G_{L'}}$ of pattern $T$ is modular, and every lift of $\Sb|_{G_{L'}}$ of pattern $T$ is modular.
\end{prop}
\begin{proof}
Fix $C_0\in\cA$. Let $m(T)$ be the number of distinct letters of weight $k$ in the image of $T$. We argue by induction on $m(T)$, for both residual representations at once. We write a pattern as a list of pairs (set of sheets, letter). In each application of Proposition~\ref{prop:swap} the lifts $B$, $A_0$, $B_0$ of the stated patterns exist by Lemma~\ref{lem:exist}.

\emph{$m(T)\le1$.} Let $T=(I\mapsto C',\ I^c\mapsto E)$ with $C\in\cA$, where $I$ may be empty (then take $C=C_0$) or everything. Let $A$ be a lift of $\Rb$ of pattern $T$ and $B$ a lift of $\Sb$ of pattern $(I\mapsto E,\ I^c\mapsto C')$. Take $A_0=f_C|_{G_{L'}}$, of constant pattern $C'$, and $B_0=\gt|_{G_{L'}}$, of constant pattern $E$; they are modular. On every sheet the two pairs are $\{C',E\}$. So $A$ and $B$ are modular. For $\Sb$: let $B$ be a lift of $\Sb$ of pattern $T$ and $A$ a lift of $\Rb$ of pattern $(I\mapsto E,\ I^c\mapsto C')$, and take again $(A_0,B_0)=(f_C|_{G_{L'}},\gt|_{G_{L'}})$. On every sheet the two pairs are $\{C',E\}$, so $A$ and $B$ are modular.

\emph{Induction step.} Let $m\ge2$ and $T=(I_0\mapsto E,\ I_1\mapsto C_1',\dots,I_m\mapsto C_m')$ with the $C_j$ distinct, $I_1,\dots,I_m$ non-empty, $I_0$ possibly empty. Let $A$ be a lift of $\Rb$ of pattern $T$. Let $B$ be a lift of $\Sb$ of pattern $(I_0\mapsto C_1',\ I_1\cup\dots\cup I_m\mapsto E)$. Let $A_0$ be a lift of $\Rb$ of pattern $(I_0\cup I_1\mapsto C_1',\ I_2\cup\dots\cup I_m\mapsto E)$ and $B_0$ a lift of $\Sb$ of pattern $(I_0\cup I_1\mapsto E,\ I_2\mapsto C_2',\dots,I_m\mapsto C_m')$. They are modular by the induction hypothesis, since their patterns have $1$ and $m-1$ letters of weight $k$. On $I_0$ the pairs are $\{E,C_1'\}$ for $(A,B)$ and $\{C_1',E\}$ for $(A_0,B_0)$; on $I_1$ both are $\{C_1',E\}$; on $I_j$, $j\ge2$, they are $\{C_j',E\}$ and $\{E,C_j'\}$. In each pair the two letters have different weights. Proposition~\ref{prop:swap} applies. For $\Sb$: let $B$ be a lift of $\Sb$ of pattern $T$, $A$ a lift of $\Rb$ of pattern $(I_0\mapsto C_1',\ I_1\cup\dots\cup I_m\mapsto E)$, $B_0$ a lift of $\Sb$ of pattern $(I_0\cup I_1\mapsto C_1',\ I_2\cup\dots\cup I_m\mapsto E)$ and $A_0$ a lift of $\Rb$ of pattern $(I_0\cup I_1\mapsto E,\ I_2\mapsto C_2',\dots,I_m\mapsto C_m')$; the pairs agree as before, and $A_0$, $B_0$ are modular by the induction hypothesis.
\end{proof}

\section{The cycles of Steinberg weights}\label{sec:St}
In \cite[Lemma 3.7(3)]{MF} the cycle of the patched module of a Serre weight is identified with a product of local cycles when the weight is not Steinberg. The proof writes the weight as a virtual combination of reductions of tame types, and tame types never contain a Steinberg weight. In this section we extend the identification to all Serre weights, using one family of wildly ramified types. The cycles of Steinberg weights can share components with cycles of other weights \cite[Theorem 1.2]{CEGS}. The argument here does not require disjointness or coefficient one.

Let $K$ be a finite extension of $\Q_p$, with uniformizer $\varpi_K$ and residue field $k$ of order $q$, $G=\GL_2(k)$, $Z$ its centre, $T\subset G$ a non-split torus, identified with $k_2^\times$ for the quadratic extension $k_2$ of $k$, and $N\colon k_2^\times\to k^\times$ the norm. Let $\kappa$ be the character of order $2$ of $k^\times$. We do not distinguish between a character of $k^\times$ with values in $\overline\F_p^\times$ and its Teichm\"uller lift.

\begin{lem}\label{lem:StGroth}
Let $\eta$ be a character of $k^\times$. In the Grothendieck group of representations of $G$ over $\overline\F_p$, the reduction of $\operatorname{Ind}_T^G(\eta\circ N)$ is the sum of the Steinberg weight $\mathrm{St}\otimes((\eta\kappa)\circ\det)$ and of non-Steinberg Serre weights.
\end{lem}
\begin{proof}
The irreducible representations of $G$ in characteristic zero are the characters $\mu\circ\det$, the representations $\mathrm{St}_\mu=\mathrm{St}\otimes\mu\circ\det$, the irreducible principal series and the cuspidal representations. The reductions of the last two kinds are reductions of tame types and contain no Steinberg weight \cite[Lemma A.4(2), proof]{CEGSb}; a character reduces to a non-Steinberg weight; and $\mathrm{St}_\mu$ reduces to the Steinberg weight $\mathrm{St}\otimes\mu\circ\det$. By Frobenius reciprocity the multiplicity of $\mathrm{St}_\mu$ in $\operatorname{Ind}_T^G\theta$ is the multiplicity of $\theta$ in $\mathrm{St}_\mu|_T$. The group $T/Z$ acts simply transitively on $\mathbf P^1(k)$, so the restriction to $T$ of $\operatorname{Ind}_B^G1=1\oplus\mathrm{St}$ is the sum of the characters of $T$ trivial on $Z$, and
\[
\mathrm{St}_\mu|_T=\bigoplus_{\theta|_Z=\mu^2,\ \theta\ne\mu\circ N}\theta .
\]
For $\theta=\eta\circ N$ the condition is $\mu^2=\eta^2$ and $\mu\circ N\ne\eta\circ N$. As $N$ is surjective this means $\mu=\eta\kappa$.
\end{proof}

\begin{lem}\label{lem:Sttype}
For every character $\eta$ of $k^\times$ there is an inertial type $\tau_\eta$ of $K$, whose tame determinant satisfies $(\det\tau_\eta)\circ\operatorname{Art}_K|_{\cO_K^\times}=\eta^2\circ\mathrm{red}$, such that the semisimplified reduction of $\sigma(\tau_\eta)$ is the inflation to $\GL_2(\cO_K)$ of the semisimplified reduction of $\operatorname{Ind}_T^G(\eta\circ N)$.
\end{lem}
\begin{proof}
Let $E$ be the unramified quadratic extension of $K$, embedded in $M_2(K)$ with $\cO_E\subset M_2(\cO_K)$, and $U^1=1+\varpi_KM_2(\cO_K)$. Choose $\beta\in\cO_E$ of trace zero whose reduction generates $k_2$ over $k$, and a non-trivial additive character $\psi_0$ of $k$. Then $1+\varpi_Kx\mapsto\psi_0(\operatorname{tr}(\bar\beta\bar x))$ is a character $\psi_\beta$ of $U^1$, because $U^1/(1+\varpi_K^2M_2(\cO_K))\cong M_2(k)$ through $1+\varpi_Kx\mapsto\bar x$; it is normalized by $E^\times$ and trivial on $1+\varpi_K\cO_K$, since $\operatorname{tr}\beta=0$. Let $\Lambda$ be the character of $J^0=\cO_E^\times U^1$ which is $\psi_\beta$ on $U^1$ and $\eta\circ N$ on the Teichm\"uller representatives $k_2^\times\subset\cO_E^\times$. Extend $\Lambda$ to $J=E^\times U^1$, writing the extension as $\widetilde\Lambda$, and put
\[
\kappa_\eta=\operatorname{Ind}_{J^0}^{\GL_2(\cO_K)}\Lambda,
\qquad
\pi_\eta=\mathrm{c\!\!-\!Ind}_{J}^{\GL_2(K)}\widetilde\Lambda^\vee.
\]
The representation $\pi_\eta$ is irreducible supercuspidal, by the unramified simple stratum construction \cite[\S\S15,19]{BH06}; the case used here is the case $r=0$ of \cite[\S5B1]{FMP}. Set $\tau_\eta=r_p(\pi_\eta)|_{I_K}$, where the unramified normalization of local Langlands does not affect inertia. These constructions are valid for an arbitrary finite extension $K/\Q_p$. Henniart's appendix to \cite{BM02} is written for an arbitrary non-archimedean local field \cite[A.1.1]{BM02}, and the pair $(J,\widetilde\Lambda)$ is the one of \cite[A.3.3--A.3.4]{BM02} for $n=1$: there $b=\varpi_K^{-1}u$ with $u$ a unit of $E$ whose residue class generates $k_2$ over $k$, the additive character is trivial on the maximal ideal and not on $\cO_K$, $J=E^\times U^1$, and the inducing representation is a character; the compact induction is irreducible supercuspidal, minimal of conductor exponent $4$. See also \cite[\S2.2]{We09}. By \cite[A.3.1]{BM02} the induction to $\GL_2(\cO_K)$ of the restriction of this character to $J^0$ is irreducible, is a type for the inertial class of the compact induction, and occurs in it with multiplicity one. Together with the dual convention of \cite[Theorem 2.1.3]{GK}, which is stated for any $K$, this gives
\[
\sigma(\tau_\eta)=\kappa_\eta.
\]
The dual in the definition of $\pi_\eta$ is required for this equality. The values of $\psi_\beta$ are $p$-th roots of unity, so the reduction of $\Lambda$ is trivial on $U^1$ and equals $\eta\circ N$ on $k_2^\times=T$. Thus the reduction of $\kappa_\eta$ is the inflation of the reduction of $\operatorname{Ind}_T^G(\eta\circ N)$. Its central character on $\cO_K^\times$ is $\eta^2\circ\mathrm{red}$. Indeed the central character of $\pi_\eta$ on units is $\eta^{-2}\circ\mathrm{red}$. Local Langlands identifies it with $(\det\tau_\eta)\circ\operatorname{Art}^{\rm geom}_K$ on units. Our arithmetic Artin map is the inverse of the geometric one, giving $\eta^2\circ\mathrm{red}$ as asserted.

\end{proof}

Now let $\bar r$ and the rings be as in \cite[\S2.3]{MF}, where the cycles $\cC_\sigma$ of the non-Steinberg weights are defined, and let $R^{\tau_\eta}$ be the potentially Barsotti--Tate lifting ring of type $\tau_\eta$ with the determinant fixed as there. For a Steinberg weight $\sigma=\mathrm{St}\otimes\mu\circ\det$ whose central character is compatible with $\det\bar r$, put, with $\eta=\mu\kappa$,
\[
\cC_\sigma:=Z(R^{\tau_\eta}/\varpi)-\sum_{\sigma'}m_{\sigma'}(\tau_\eta)\,\cC_{\sigma'},
\]
the sum being over the non-Steinberg weights and $m_{\sigma'}(\tau_\eta)$ the multiplicity of $\sigma'$ in the reduction of $\sigma(\tau_\eta)$.

\begin{prop}\label{prop:Stcycle}
In a global setting, for every Serre weight $\sigma=\bigotimes_v\sigma_v$ with the central character fixed in \cite[\S3]{MF},
\[
Z(M_\infty^\sigma)=2\prod_{v\mid p}\cC_{\sigma_v}\times\cD_\Sigma ,
\]
where $\cD_\Sigma$ is the cycle of \cite[Lemma 3.7]{MF} coming from the places of $\Sigma$ and the patching variables.
\end{prop}
\begin{proof}
For families of non-Steinberg weights this is \cite[Lemma 3.7(3)]{MF}. The proof given there applies to all weights once the types $\tau_\eta$ are allowed. First, \cite[Proposition 3.4(3),(4)]{MF} hold for a family $\tau=(\tau_v)$ in which each $\tau_v$ is a tame type or one of the types $\tau_\eta$: the category on which the patching functor is defined contains the lattices in $\bigotimes_v\sigma(\tau_v)$; the proof of (3) uses local-global compatibility and the characterization of $\sigma(\tau_v)$ by the inertial local Langlands correspondence, for any inertial type; and (4) rests on \cite[Lemma 5.3.6]{GK}, which is proved for arbitrary inertial types. Hence the identity in the proof of \cite[Lemma 3.7(3)]{MF},
\[
\prod_{v\mid p}Z(R_v^{\tau_v}/\varpi)\times\cD_\Sigma=\tfrac12\sum_{\sigma'}\ \prod_{v\mid p}m_{\sigma'_v}(\tau_v)\,Z(M_\infty^{\sigma'}),
\]
holds for these families. Second, by Lemmas~\ref{lem:StGroth} and~\ref{lem:Sttype} and by the definition of the cycles, every Serre weight $\sigma_v$ is a combination $\sum_\tau n_\tau(\sigma_v)\overline{\sigma(\tau)}$ with integer coefficients, over the tame types and the types $\tau_\eta$, such that $\cC_{\sigma_v}=\sum_\tau n_\tau(\sigma_v)Z(R_v^\tau/\varpi)$. Multiplying the identity by $\prod_vn_{\tau_v}(\sigma_v)$ and summing over the families gives the proposition, as in \cite[Lemma 3.7(3)]{MF}. The integral spanning assertion used here also follows directly from the proof of \cite[Lemma~4.5.2 and Corollary~4.5.3]{GK}, which treats every finite $K/\Q_p$ and separates the summands with a fixed tame determinant.
\end{proof}

This is the argument of \cite[Theorem 5.2]{CEGS} with wildly ramified types, as indicated in \cite[Remark 5.3]{CEGS}. A related patching comparison, under its stated globalization hypotheses, is \cite[Theorem 5.5.2]{EG}.

\section{Proof of Theorem \ref{thm:main}}\label{sec:proof}
A \emph{word} is a family $(C^{(1)},\dots,C^{(d)})\in\cA^d$; it defines a family of components over $L$, with $C^{(i)}_v$ at the place of the sheet $i$ above $v$.

\subsection*{Every word is automorphic over $L$} If $\cA$ is empty there is no word, since $d\ge1$. Otherwise fix a word. By (I2) over $L$, there is a lift $\rho_T$ of $\Rb|_{G_L}$ with determinant $\det\rho|_{G_L}$, unramified outside $p$ and the places above $\Sigma$, whose points above $p$, after the twist $\nu_{\mathrm{MF}}$, lie on the components of the word, and whose points above $\Sigma$ lie on the prescribed Steinberg rings. Nonzero monodromy at these latter places will follow after modularity is proved. Its restriction to $G_{L'}$ is a lift of the pattern $i\mapsto (C^{(i)})'$, so it is modular over $L'$ by Proposition~\ref{prop:mixed}. As $L'/L$ is soluble and the restriction is irreducible, $\rho_T$ is modular over $L$. At a place $w$ above $\Sigma$, $\rho_T$ corresponds to a point of the Steinberg ring; this ring also has unramified points, with Frobenius eigenvalues of ratio the cardinality of the residue field, but for a modular representation these are excluded by local-global compatibility and the Ramanujan conjecture \cite{Bla06}, as in \cite[Lemma 3.3]{MF}. So inertia at $w$ acts unipotently and non-trivially, and $\rho_T$ has the properties $(\dagger)$ over $L$. By \cite[Lemmas 3.3 and 3.6(1)]{MF} the word is in $\cA_L$.

\subsection*{Counting} For $v\mid p$ let
\[
s_v=\sum_{\sigma_v\in \mathrm{JH}(W_{k,v})}n_{k,v}(\sigma_v)\,\cC_{\sigma_v},
\]
where $n_{k,v}(\sigma_v)$ is the multiplicity of $\sigma_v$ in $W_{k,v}$ and $\cC_{\sigma_v}$ is the cycle of \cite[\S2.3]{MF} if $\sigma_v$ is not Steinberg, and the cycle of Section~\ref{sec:St} if it is. The multiplicity of $\bigotimes_v\sigma_v$ is $\prod_v n_{k,v}(\sigma_v)$. Summing Proposition~\ref{prop:Stcycle} with these multiplicities and using \cite[Lemma~3.7(1)]{MF} gives the following identities. In each setting, multiplication by its own nonzero factor $2\cD_\Sigma$ is injective on the subgroup spanned by the product cycles, by \cite[Lemma~2.6]{MF}.
\[
\sum_{C\in\cA}\ \prod_{v\mid p}Z(\overline C_v)=\prod_{v\mid p}s_v,\qquad
\sum_{T\in\cA_L}\ \prod_{u\mid p}Z(\overline T_u)=\prod_{u\mid p}s_{v(u)} .
\]
The local data over $L$ are those over $F$, and the right-hand side of the second identity is the product over the sheets of $\prod_vs_v$. The product over the sheets of the first identity therefore shows that the sum of $\prod_uZ(\overline T_u)$ over the $|\cA|^d$ words equals the right-hand side of the second identity. Distinct words are distinct families, and all words are in $\cA_L$. The cycles $\cC_{\sigma_v}$ are defined by the local rings at $v$ alone, which is all that is used. Every term of the left-hand side is a non-zero effective cycle, so there is no other term, and $\cA_L$ is exactly the set of words.

\subsection*{Conclusion} Write $C_{\rho,v}$ for the component of the ring of \cite{MF} containing $(\rho\otimes\nu_{\mathrm{MF}})|_{G_{F_v}}$. $\rho$ is modular over $L$ and has the properties $(\dagger)$ by \S\ref{ss:red}, so the family $(C_{\rho,v(u)})_u$ of its components is in $\cA_L$ \cite[Lemmas 3.3 and 3.6(1)]{MF}. Hence it is a word, and $(C_{\rho,v})_v\in\cA$. By \cite[Lemma 3.6(2)]{MF}, $\rho$ comes from an eigenform contributing to $M_k$, and the theorem follows from \S\ref{ss:red}. \qed

\section{The general case}\label{sec:general}
We indicate the changes for Theorem~\ref{thm:general}. The framework is that of \cite[\S6]{MF}: the global settings, the lattices $L_{\lambda,\tau}$, the modules $M_{\lambda,\tau}$, the properties $(\dagger_{\lambda,\tau})$ and the set $\cA$ are as there, and \cite[Lemmas 3.6 and 3.7]{MF} hold with the multiplicity $n_{\lambda,\tau}$ in place of $n_k$.

\emph{Weights.} Let $w$ be the sum of the two Hodge--Tate weights of $\rho$ at an embedding; it does not depend on the embedding, because $\det\rho$ is a character of the Galois group of a totally real field. Let $\mathcal K$ be the set of the integers $k_\kappa$, where $k_\kappa-1$ is the difference of the two Hodge--Tate weights at the embedding $\kappa$; they all have the parity of $w+1$. Take an auxiliary form for $(\Rb,w,\mathcal K)$. At every embedding the four Hodge--Tate weights of $\rho\otimes\gt$ are distinct, because $k_S\neq k_\kappa$.

\emph{Fields.} In \S\ref{sec:setup} we require in addition that, for $w\mid p$, $\tau_v|_{I_{F'_w}}$ be trivial. Then the restrictions to $G_{K'_v}$ of the points of a component of the potentially crystalline ring of type $(\lambda_v,\tau_v)$ are crystalline, and the rings $R^k_v$ and $R^{(4)}_v$ are replaced by the crystalline rings over $K'_v$ with the Hodge--Tate weights of $\rho$ and of $\rho\otimes\gt$.

\emph{Mixed lifts.} Section~\ref{sec:mixed} is unchanged.

\emph{Words.} In \S\ref{sec:proof} the lift $\rho_T$ over $L$ is taken, by (I2), on the components of the potentially crystalline rings of types $(\lambda_v,\tau_v)$ given by the word. To apply (I2), first undo the common patching twist. For a finite Galois extension killing $\tau_v$ on inertia, the inertial type is locally constant on the generic fibre with semistability over that extension, as can be seen from the traces of its finite inertia quotient. Inside the fixed type, the crystalline quotient and the semistable ring have the same generic dimension \cite[Theorems~3.3.4 and~3.3.8]{Kis08}. Each crystalline component therefore gives a component allowed by (I2). The counting uses the cycles $s_v$ defined with the multiplicities of the Serre weights in $\overline{\sigma(\lambda_v,\tau_v)}$, and Proposition~\ref{prop:Stcycle}. \qed

\section{Potentially semistable representations}\label{sec:ss}
For a regular Hodge type $\lambda_v$ and an inertial type $\tau_v$ at $v\mid p$, let $\sigma^{\mathrm{ss}}(\tau_v)$ be the representation of $\GL_2(\cO_{F_v})$ attached to $\tau_v$ by Henniart and Kisin for potentially semistable representations: it is $\sigma(\tau_v)$ if $\tau_v$ is not scalar, and $\sigma(\tau_v)\otimes\widetilde{\mathrm{St}}$ if $\tau_v$ is scalar, where $\widetilde{\mathrm{St}}$ is the inflation of the Steinberg representation of $\GL_2(k_v)$. Let $\sigma^{\mathrm{ss}}(\lambda_v,\tau_v)$ be its tensor product with the algebraic representation of weight $\lambda_v$, and $R^{\lambda_v,\tau_v,\mathrm{ss}}_v$ the lifting ring of potentially semistable lifts of Hodge type $\lambda_v$ and inertial type $\tau_v$ \cite{Kis08}. For scalar $\tau_v$, this type occurs once in both generic cases: after untwisting, the first principal congruence invariants of an unramified principal series are $1\oplus\widetilde{\mathrm{St}}$, and those of an unramified Steinberg representation are $\widetilde{\mathrm{St}}$. For nonscalar $\tau_v$, the monodromy is necessarily zero. Its points include the potentially crystalline ones; by \cite[Theorem 3.3.4]{Kis08} its generic fibre is equidimensional, but it is not regular in general.

\begin{thm}\label{thm:ss}
Let $p\ge7$ be a prime and $F$ a totally real field. Let $\rho\colon G_F\to\GL_2(\Qbar_p)$ be a continuous representation, unramified outside finitely many places, with reduction $\Rb$. Assume:
\begin{enumerate}[label=\textup{(\arabic*)}]
\item for every $v\mid p$, $\rho|_{G_{F_v}}$ is potentially semistable with distinct Hodge--Tate weights at every embedding; let $\lambda_v$ be its Hodge type and $\tau_v$ the restriction to inertia of its Weil--Deligne representation;
\item $\rho$ is potentially modular;
\item $\Rb$ is modular, and $\Rb|_{G_{F(\zeta_p)}}$ is absolutely irreducible with non-dihedral projective image.
\end{enumerate}
Then $\rho$ is modular: it comes from a Hilbert modular eigenform over $F$ of the weight given by the $\lambda_v$, whose local component $\pi_v$ at each $v\mid p$ contains $\sigma^{\mathrm{ss}}(\tau_v)^\vee$.
\end{thm}

\begin{cor}\label{cor:ss}
Let $p\ge7$, $k\ge2$, and $\rho$ as in Theorem~\ref{thm:main}, with hypothesis (1) replaced by: for every $v\mid p$, $\rho|_{G_{F_v}}$ is semistable with Hodge--Tate weights $0$ and $k-1$ at every embedding. Then $\rho$ comes from a Hilbert modular eigenform over $F$ of parallel weight $k$ whose level is prime to $p$ times an Iwahori level at the places $v\mid p$ at which $\rho$ is not crystalline; at those places $\pi_v$ is an unramified twist of the Steinberg representation.
\end{cor}
\begin{proof}
This is Theorem~\ref{thm:ss} with trivial inertial types. At a place where $\rho$ is not crystalline, the monodromy of its Weil--Deligne representation is non-zero, so by local-global compatibility \cite{Sai09,Ski09} $\pi_v$ is an unramified twist of Steinberg, which has conductor exponent one and has Iwahori fixed vectors. At a crystalline place both the inertial type and the monodromy are trivial. Local-global compatibility and local Langlands then give a spherical local component, as required.
\end{proof}

The proof of Theorem~\ref{thm:ss} is that of Theorem~\ref{thm:general}, with $\sigma^{\mathrm{ss}}(\lambda_v,\tau_v)$ in place of $\sigma(\lambda_v,\tau_v)$ and semistable rings in place of crystalline ones. Smoothness of the generic fibres of the local rings was used in three places, and we explain how each is replaced. We first record a lemma.

\begin{lem}\label{lem:geomint}
Let $R$ be a complete local noetherian $\cO$-algebra with residue field $\F$, and $\mathfrak p_1,\dots,\mathfrak p_r$ the minimal primes of $R$ with $R/\mathfrak p_i$ flat over $\cO$. There is a finite extension $E'/E$ such that for every minimal prime $\mathfrak q$ of $R\otimes_\cO\cO_{E'}$ which is flat over $\cO_{E'}$, the ring $(R\otimes_\cO\cO_{E'})/\mathfrak q$ is geometrically integral in the sense of \cite[\S3]{BLGHT}.
\end{lem}
\begin{proof}
It suffices to treat one $A=R/\mathfrak p_i$. Its normalization $B$ in $K_A=\operatorname{Frac}(A)$ is finite, because $A$ is excellent. Since $A$ is complete local and $B$ is a domain, $B$ is complete local, with residue field $\widetilde\F$ finite over $\F$. Choose a height one prime $\mathfrak P$ of $B$ containing the coefficient uniformizer $\varpi_E$. The normal local ring $B_{\mathfrak P}$ is a discrete valuation ring.

Let $E_A$ be the algebraic closure of $E$ in $K_A$, and let $E''/E$ be any finite subextension. Its ring of integers is contained in $B$, since its elements are integral over $\cO\subset A$. If $\varpi''$ is a uniformizer and $e=e(E''/E)$, the equality $\varpi_E=u(\varpi'')^e$, with $u$ a unit of $\cO_{E''}$, gives
\[
v_{\mathfrak P}(\varpi_E)=e\,v_{\mathfrak P}(\varpi'').
\]
Both $u$ and $u^{-1}$ lie in $B$, so $v_{\mathfrak P}(u)=0$, and $v_{\mathfrak P}(\varpi'')$ is a positive integer. Thus $e$ is bounded by $v_{\mathfrak P}(\varpi_E)$. The inclusion $\cO_{E''}\subset B$ is local, so the residue degree is bounded by $[\widetilde\F:\F]$. The degrees of all finite subextensions are therefore bounded, and $E_A/E$ is finite.

Choose $E'$ containing the Galois closure of $E_A/E$. In characteristic zero, $K_A/E_A$ is regular, since $E_A$ is algebraically closed in $K_A$. Consequently $K_A\otimes_EE'$ is a finite product of fields regular over $E'$. Each flat minimal prime quotient of $A\otimes_\cO\cO_{E'}$ embeds in one of these fields. Tensoring with any further finite coefficient extension preserves that embedding and gives a domain, because the corresponding tensor product of fields is a field. Hence these quotients are geometrically integral. Take a common $E'$ for the finitely many $A$.

\end{proof}

Since only finitely many local rings occur, we enlarge $E$ once and for all so that Lemma~\ref{lem:geomint} holds for all of them. This is harmless, and was already implicit in the use of \cite[Lemma 3.3]{BLGHT} in \cite{MF} and in Lemma~\ref{lem:tensor}.

\emph{1. Products of components.} In Lemma~\ref{lem:tensor} and in \cite[Lemma 3.7(1)]{MF}, the completed tensor products of quotients of local rings by minimal primes are integral, and their minimal primes are the sums of the minimal primes of the factors. By \cite[Lemma 3.3(3),(5)]{BLGHT} this holds as soon as these quotients are $\cO$-flat and geometrically integral, which is Lemma~\ref{lem:geomint}. Smoothness is not needed.

\emph{2. The patched module.} If the patched module is zero there is nothing to prove. Otherwise, for lattices in $\bigotimes_v\sigma^{\mathrm{ss}}(\lambda_v,\tau_v)$, the module $M_\infty$ is finite free over the Taylor--Wiles power series ring $S_\infty$. The dimension equality with $R_\infty^{\mathrm{ss}}$ makes $M_\infty$ maximal Cohen--Macaulay over $R_\infty^{\mathrm{ss}}$. Thus its support is a union of irreducible components. Let $R_M=R_\infty^{\mathrm{ss}}/\operatorname{Ann}(M_\infty)$. Its faithful action embeds it in $\operatorname{End}_{S_\infty}(M_\infty)$, so $R_M$ is finite over $S_\infty$. This algebra is torsion free over $S_\infty$, so each of its minimal primes contracts to zero. Thus each supported component maps surjectively to $\operatorname{Spec}S_\infty$ and meets its characteristic zero augmentation point. The resulting point comes from an eigenform, as in \cite[Lemma B.5.1 and its proof]{GK}.

At such a point, the local component of the cuspidal $\GL_2$ representation is generic. Local compatibility and \cite[Lemma 1.1.3 and Theorem 1.2.7]{All16} show that the corresponding potentially semistable point above $p$ is formally smooth. This also holds with fixed determinant: the unrestricted determinant adds a formally smooth unramified parameter, split locally by the square root twist because $p>2$. At the other places the smoothness argument is that of \cite{MF}. Hence the patched local ring at this point is regular, and Auslander--Buchsbaum makes $M_\infty$ free there. Its reduction modulo the augmentation ideal is still faithful over the corresponding quotient ring; the reduced characteristic zero Hecke algebra and multiplicity two calculation of \cite[Proposition 3.4]{MF} then give rank $2$. Every supported component therefore has generic rank $2$.

This is enough for the identity of cycles in \cite[Lemma 3.7(1)]{MF}: \cite[Proposition 2.2.13]{EG} uses $\varpi$-torsion freeness and the lengths at the generic points of the support. It does not require smoothness of the entire generic fibre.

\emph{3. Potentially connected representations.} In Definition~\ref{def:potconn} we replace \emph{crystalline} by \emph{semistable} in (ii) and the crystalline ring by the semistable ring in (iii). After restriction to the comparison field and identification of the residual representations, a common component of this semistable ring gives the relation in \cite[Definition~3.27]{ADP}. That definition compares two lifts of the same residual representation over the field where the lifting theorem is applied. Remark~\ref{rem:potconn}(b) no longer holds as stated: a point of the generic fibre may lie on two components, so the relation need not be transitive. Transitivity is never used: Proposition~\ref{prop:swap} and (I10) only require that two representations lie on a common component. Parts (a), (c) and (d) of the remark hold, with (I10) in place of (I4) in (d).

The remaining changes are the following.

\emph{Fields.} $F'$ is chosen so that $\tau_v|_{I_{F'_w}}=1$ for $w\mid p$. The restrictions to $G_{K'_v}$ of the points of a component of $R^{\lambda_v,\tau_v,\mathrm{ss}}_v$ are semistable and form an irreducible image. Choose a component of the semistable ring over $K'_v$ containing this image. For an automorphic component the image contains the local point of $f_C|_{G_{L'}}$, which is automorphic by soluble base change. Its Weil--Deligne representation is therefore generic, so its containing component is unique by \cite[Theorem 1.2.7]{All16}. The rings $R^k_v$ and $R^{(4)}_v$ are replaced by the semistable rings over $K'_v$ with the Hodge--Tate weights of $\rho$ and of $\rho\otimes\gt$. The ring $R^g_v$ is unchanged.

\emph{Lemma~\ref{lem:tensor}.} A tensor product of semistable representations is semistable, with monodromy $N\otimes1+1\otimes N$. In this analogue of the lemma we assert existence of a common component, without uniqueness. The product of the two source components is integral, so choose a target component containing its irreducible image. The conjugation argument uses only density of the characteristic zero points and stability of the Hodge type and semistability under conjugation. It preserves each component, and therefore puts the products in the opposite order on the same chosen component.

\emph{Mixed lifts.} In Lemma~\ref{lem:exist}, (I2) allows potentially semistable components. For (I3), discard the fixed determinant by adjoining the unramified determinant parameter. At a split CM place the polarized lifting ring identifies with the usual lifting ring \cite[Remark~2.1.11]{CEG}. The image of a fixed determinant component thus lies on the component obtained by adjoining that parameter, which contains the same automorphic witness. Lifts of $\Sb$ are now taken on letters that may be semistable and not crystalline; they are still potentially modular and Lie-irreducible. Put $m=\chi_R\chi_S$ and use the character $\theta$ of Proposition~\ref{prop:swap}. Its construction also applies here, since the determinant of a semistable representation is a semistable character, hence crystalline. For the twist of both tensors the identity on $G_{E_0}$ is
\[
m\,\theta\theta^c=\eps^{-3}.
\]
The raw multiplier on $G_{L'}$ is exactly $\eps^{-3}$. Before twisting it is $m\delta_{E_0/L'}$, from the symmetric tensor pairing. Twisting multiplies it by $\theta\circ\operatorname{Ver}_{E_0/L'}$, where $\operatorname{Ver}_{E_0/L'}$ is transfer. Both sides agree on $G_{E_0}$ by the displayed identity, and on complex conjugation they both have value $-1$, since transfer sends a complex conjugation to its square. Thus the twisted target and seed have the same polarized reduction, and the polarizing Hecke character of the twisted seed can be taken to be trivial. Twisting by the same crystalline character preserves the common semistable component above $p$ after shifting the Hodge type. Adequacy and the common component condition at every other place follow as in Proposition~\ref{prop:swap}, including places ramified in $\theta$. Apply (I10) to the twisted pair, remove the Hecke character corresponding to $\theta$, and descend the resulting untwisted automorphic tensor as in that proposition.

\emph{Proof.} In \S\ref{sec:proof} the lifts $\rho_T$ are taken on components of the potentially semistable rings, the cycles $s_v$ are defined with the multiplicities of the Serre weights in $\overline{\sigma^{\mathrm{ss}}(\lambda_v,\tau_v)}$, and Proposition~\ref{prop:Stcycle} gives the cycles of all of them. Even for trivial $\tau_v$ and parallel weight, $\overline{\sigma^{\mathrm{ss}}(\lambda_v,\tau_v)}=W_{k,v}\otimes\overline{\mathrm{St}}$ can contain a Steinberg weight or a constituent of multiplicity greater than one. For any word, the constructed lift $\rho_T$ is modular over $L$ by Proposition~\ref{prop:mixed} and soluble descent. Its local points above $p$ are therefore smooth by \cite[Theorem~1.2.7]{All16} and local-global compatibility. At each such point the component prescribed by the word is the unique component through the point. The associated Hecke eigenform consequently puts the whole word in $\cA_L$. The counting argument now gives $\cA_L=\cA^d$. Finally $\rho|_{G_L}$ is modular and every $v\mid p$ splits completely in $L$. Thus the local point of $\rho$ at $v$ is smooth on the ring over $F_v$, and its unique component occurs at every place above $v$. The corresponding constant family over $L$ is a word, so the family over $F$ is in $\cA$. The point defined by $\rho$ lies in the support of $M_\infty$ above the augmentation ideal, and specialization gives an eigenform over $F$ as in the proof of \cite[Lemma~3.6(2)]{MF}. \qed

\section{Remarks}
\begin{rem}[Relation with \cite{MF}]
For $p\ge7$ and non-dihedral projective image of $\Rb|_{G_{F(\zeta_p)}}$, Theorems~\ref{thm:main} and \ref{thm:general} contain \cite[Theorems 1.2 and 1.7]{MF} and the variants under a uniqueness hypothesis \cite[Propositions 4.4 and 6.1]{MF}. The proof uses \cite{MF} through the reduction to a global setting, the identity of cycles, and the two statements of \cite[Lemma 3.6]{MF}. The multiplicity-free argument of \cite{MF} is the case in which no mixed lift is needed.
\end{rem}

\begin{rem}[Automorphic and potentially automorphic components]
The identity $\cA_L=\cA^d$ obtained in \S\ref{sec:proof} says that, in the situation of the theorem, the set of automorphic components over the chosen extension $L/F$ of Section~\ref{sec:setup} is the product, over its sheets, of the sets over $F$. In particular, under the hypotheses of the theorem, a family of local components which contains the local representations of a potentially modular lift with the properties $(\dagger)$ is automorphic.
\end{rem}

\begin{rem}[What is not covered]
The hypotheses require $p\ge7$ and absolute irreducibility with non-dihedral projective image after restriction to $G_{F(\zeta_p)}$.
\end{rem}

\subsection*{Acknowledgments}
The author used Claude, an AI model developed by Anthropic, as an assistant during the preparation of this paper: for discussing and checking the arguments, for searching and verifying the literature, for adversarial reviews of successive drafts, and for help with writing and typesetting. All mathematical statements and proofs were checked by the author, who takes full responsibility for the content of the paper.

\end{document}